\documentclass[11pt]{article}

\usepackage{graphicx}
\usepackage{latexsym,amsmath,amsfonts,amscd, amsthm, dsfont}
\usepackage{mathrsfs}
\usepackage{bm,color}
\usepackage{epsfig,verbatim,epstopdf,graphics}
\usepackage{subfigure}
\usepackage{changebar}
\usepackage{multirow}

\usepackage{algorithmic}

\usepackage{yhmath}
\usepackage{booktabs} 
\usepackage{tikz}
\usepackage{verbatim}
\usetikzlibrary{arrows,backgrounds,snakes,shapes}

\usepackage{xcolor}
\usepackage{xpatch}

\makeatletter
\newcommand*{\rom}[1]{\expandafter\@slowromancap\romannumeral #1@}
\makeatother

\makeatletter 
\let\myorg@bibitem\bibitem
\def\bibitem#1#2\par{%
	\@ifundefined{bibitem@#1}{%
		\myorg@bibitem{#1}#2\par
	}{%
		\begingroup
		\color{\csname bibitem@#1\endcsname}%
		\myorg@bibitem{#1}#2\par
		\endgroup
	}%
}
\newcommand*{\bibitem@renhigh}{blue}    
\newcommand*{\bibitem@HWENOAP}{blue}
\newcommand*{\bibitem@PaillereDeconinckRoe}{blue}
\newcommand{\bibitem@ShuOsherJCPO}{blue}
\newcommand{\bibitem@ShuOsherJCPT}{blue}
\newcommand*{\bibitem@mythesis}{orange}  
\makeatother 

\numberwithin{equation}{section}

\graphicspath{{./}{./figure/}}
\allowdisplaybreaks

\newtheoremstyle{plainNoItalics}{}{}{\normalfont}{}{\bfseries}{.}{ }{}

\theoremstyle{plain}
\newtheorem{thm}{Theorem}[section]

\theoremstyle{plainNoItalics}

\newtheorem{rem}[thm]{Remark}

\newtheorem{exa}[thm]{Example}

\newcommand{\be}{\begin{eqnarray}}
\newcommand{\ee}{\end{eqnarray}}
\newcommand{\beno}{\begin{eqnarray*}}
\newcommand{\eeno}{\end{eqnarray*}}

\makeatletter

\newcommand{\Rmnum}[1]{\expandafter\@slowromancap\romannumeral #1@}
\makeatother

\begin{document}

\baselineskip=1.8pc


\begin{center}
{\Large \bf High Order Residual Distribution Conservative Finite Difference  HWENO Scheme for Steady State Problems}
\end{center}

\vspace{.2in}
\centerline
{
Jianfang Lin \footnote{School of Mathematical Sciences, Hangzhou, Zhejiang University, Zhejiang 310058, P.R. China. {\tt jianfang.lin@zju.edu.cn}},
Yupeng Ren \footnote{School of Mathematical Sciences, Xiamen University, Xiamen, Fujian 361005, P.R. China. {\tt ypren@stu.xmu.edu.cn}},
 R\'emi Abgrall \footnote{Institute of Mathematics, University of Zurich, Zurich 8057, Switzerland. {\tt remi.abgrall@math.uzh.ch}},
Jianxian Qiu\footnote{School of Mathematical Sciences and Fujian Provincial Key Laboratory of Mathematical Modeling and High-Performance Scientific Computing, Xiamen University, Xiamen, Fujian 361005, P.R. China. {\tt jxqiu@xmu.edu.cn}}
}

\bigskip
\noindent
{\bf Abstract.}
In this paper, we develop a high order residual distribution (RD) method for solving steady state conservation laws in a novel Hermite weighted essentially non-oscillatory (HWENO) framework recently developed in \cite{renhigh}. In particular, {\color{orange}we design a high order HWENO integration for the integrals of source term and fluxes based on the point value of the solution and its spatial derivatives,} and the principles of residual distribution schemes are adapted to obtain steady state solutions. {\color{orange}Two advantages of the novel HWENO framework have been shown in \cite{renhigh}: first, compared with the traditional HWENO framework, the proposed method does not need to introduce additional auxiliary equations to update the derivatives of the unknown variable, and just compute them from the current point value of the solution and its old spatial derivatives, which saves the computational storage and CPU time, and thereby improve the computational efficiency of the traditional HWENO framework.} Second, compared with the traditional WENO method, reconstruction stencil of the HWENO methods becomes more compact, their boundary treatment is simpler, and the numerical errors are smaller at the same grid. Thus, it is also a compact scheme when we design the higher order accuracy, compared with that in \cite{Chou.Shu_JCP2006} Chou and Shu proposed. {\color{orange}Extensive numerical experiments for one- and two-dimensional scalar and systems problems confirm the high order accuracy and good quality of our scheme.}


\vfill

{\bf Key Words:} High order accuracy; Residual distribution; HWENO scheme; Conservation laws; Steady state.
\newpage

\section{Introduction}
In this paper, we propose a new type of the residual distribution (RD) conservative finite difference Hermite weighted essentially non-oscillatory (HWENO) method for solving the following steady state hyperbolic conservation laws
\begin{equation}\label{Eq:ssp}
\nabla\cdot F(u) = 0,
\end{equation}
{\color{blue}where hyperbolicity means that $\frac{\partial F(u)}{\partial u}$ is diagonalizable with real eigenvalues.} The RD scheme has received considerable attention and been successfully used to solve steady state problems in \cite{PaillereDeconinckRoe,  Abgrall_JCP2001, Abgrall.Marpeau_JSC2007}.

The RD schemes are composed of two parts: residual (or fluctuation) evaluation and residual distribution. {\color{orange} A brief framework of the RD scheme for a two-dimensional steady state  problem \eqref{Eq:ssp} is introduced as follows:} given a general triangular or quadrilateral mesh $\mathscr{T}_{h}$, nodes $\left\{M_{i}\right\}_{i=1, \cdots, n_{s}}$ of $\mathscr{T}_{h}$, and $T$ is a generic element. On each element $T$, {\color{orange}we should define a total residual $\Phi^{T}$, and $\Phi^{T}_{i}$ which is the amount of $\Phi^{T}$ associated with the vertex $M_{i}$,} such that a conservation property is satisfied
\begin{equation}
\Phi^{T} = \int_{T}\! \nabla\cdot F^{h}(u_{h}) \, \mathrm{d}x ,~~\sum\limits_{i, M_{i}\in T}\Phi^{T}_{i} = \Phi^{T},
\end{equation}
{\color{blue}thus the residual distribution scheme for the two-dimensional steady state problem \eqref{Eq:ssp}  is given as
\begin{equation}
u^{n+1}_{i}=u^{n}_{i}-\frac{\Delta t_{n}}{|C_{i}|}\sum\limits_{T\in M_{i}, i}\Phi^{T}_{i},
\end{equation}
where $|C_{i}|$ is the area of the control volume associated with the vertex $M_{i}$.} The accuracy can be obtained at steady state when vanishing cell residuals. The class of RD schemes, or fluctuation splitting schemes were pioneered and developed by Roe, Sidikover, Deconinck, Struijs and their collaborators {\color{blue}\cite{Struijs.Deconinck.Roe_CFD1991, Roe.Sidilkover_JNA1992, Deconinck.Struijs.Bourgeois.Roe_CFD1993, Abgrall.Mer.Nkonga_2002, Abgrall.Roe_JSC2003, Csik.Deconinck_IJNMF2002, Abgrall.Mezine_JCP2003, Abgrall.Mezine_JCP2004, Abgrall.Marpeau_JSC2007}}. In the past decades of the development, the RD scheme has demonstrated its robustness in many numerical experiments and does not have the restriction on the regularity of the mesh. {\color{blue}The Lax-Wendroff theorem in \cite{Abgrall.Mer.Nkonga_2002} has proved that the numerical solution of the RD scheme is convergent to the weak solution, if the flux function satisfies the Lipschitz continuity. And the stability of the RD scheme can be obtained by the maximum principle, see \cite{Abgrall_JCP2001, Abgrall.Mezine_JCP2004}. The accuracy of the scheme is reached at steady state when residues vanish, if the residual property in \cite{Abgrall_JCP2001} holds}. Above the work of the RD scheme, it is at most second order accuracy. {\color{red} Later on, Abgrall and Roe \cite{Abgrall.Roe_JSC2003}  extended the RD scheme to a higher order on triangle meshes.} Also, {\color{orange} Abgrall and Meapeau \cite{Abgrall.Marpeau_JSC2007} considered the construct of the second order RD scheme on quadrilateral meshes.} Besides for the steady sate problems, the RD scheme is also applied for solving unsteady state problems, see \cite{Abgrall.Mezine_JCP2003, Csik.Deconinck_IJNMF2002}. {\color{red} Abgrall and his collaborators in \cite{Abgrall.Bacigaluppi.Tokareva_2016, Abgrall_JSC2017}  extended the RD scheme to the multi-dimensional systems.}

Due to the fact that the steady state problem \eqref{Eq:ssp} would exhibit hyperbolic behavior, for instance, shock and other discontinuities, it's necessary to develop a numerical scheme with capable of capturing these traits. Weighted essentially non-oscillatory (WENO) scheme {\color{blue}\cite{ShuOsherJCPO, ShuOsherJCPT, Jiang.Shu_JCP1996}} has been widely studied to solve hyperbolic conservation laws with good properties of high order accuracy in smooth regions and non-oscillatory near discontinuity. {\color{red} In particular, Chou and Shu in \cite{Chou.Shu_JCP2006, Chou.Shu_JCP2007}  proposed a new RD scheme combined with WENO scheme.} {\color{orange} Recently, high order Hermite WENO (HWENO) methods \cite{qiushuhwenohj2005,lqfdhweno1} have gained much attention in solving hyperbolic conservation laws.} Both the traditional WENO and HWENO methods can achieve the high order accuracy and preserve the essentially non-oscillatory property,
but the HWENO scheme uses the Hermite interpolation in reconstructing polynomials, {\color{orange}which involves both the unknown variable and its first order spatial derivative or first moment.} Thus, the HWENO reconstruction stencil becomes more compact and their boundary treatment is much simpler, although more storage and some additional work are needed to evaluate the spatial derivatives.
The HWENO scheme was first introduced in the construction of a limiter for the DG method \cite{hwenolim1,hwenolim2} due to its compact stencil. It was first used to solve the time-dependent Hamilton-Jacobi equation in \cite{qiushuhwenohj2005}, {\color{orange}and the numerical results show that the HWENO scheme has smaller errors than the traditional WENO method on the same mesh and the same order of accuracy.} Since then, {\color{orange}many HWENO schemes have been developed to solve hyperbolic conservation laws on structured and triangular meshes, see \cite{lqfdhweno1,taohwenosta,7thhweno,zhaohyhweno,zhu2018new}.} In addition, it is observed in \cite{HWENOAP} that the finite volume HWENO scheme enjoys the asymptotic preserving property, when applied to the steady-state discrete ordinates ($S_{N}$) transport equation.

In this paper, built upon the high order RD schemes and HWENO reconstructions, {\color{orange}we design a sixth order RD finite difference conservative HWENO scheme for solving steady state hyperbolic conservation laws.} It is worth mentioning that we are no longer using the traditional HWENO framework, namely, {\color{orange}using one equation to update original variables $u$ and several additional auxiliary equations to update their derivatives.} Instead, {\color{orange}we develop a novel HWENO framework, motivated by the work \cite{renhigh}, which only  uses one equation to update original variable and  in which the derivative of $u$  is obtained by applying the HWENO reconstruction on the updated values of $u$ and the previous values of the derivatives. The novel HWENO framework not only inherits the advantages of the traditional HWENO, but also  saves the computational storage and CPU costs, which improves the computational efficiency of the traditional HWENO scheme. Moreover,  for two-dimensional steady state hyperbolic conservation laws, it is  still not clear and difficult to distribute the residuals for the auxiliary equations under the traditional HWENO framework and will be explored in future.}

The paper is organized as follows. In Section $2$ and $3$, {\color{orange}we formulate the residual evaluation and the residual distribution procedures for one- and two-dimensional problems, respectively. The performances of the proposed method are demonstrated in Section $4$,} through extensive numerical tests on several benchmark problems for steady state simulations. Finally, concluding remarks are given in Section $5$.

\section{High order RD conservative finite difference HWENO scheme in one dimension}
In this section, {\color{orange}we develop a high order RD conservative finite difference HWENO scheme for steady state hyperbolic conservation laws in one dimension.} In the first subsection, we define the total residual within each interval through the integral form, and then describe the distribution of total residual within each interval, complying with the principles of upwind scheme and the residual property. In the second subsection, we generalize the scheme to one-dimensional systems, based on a local characteristic field decomposition and using the same principles as in the scalar case to distribute the total residual within each interval in the characteristic fields.

\subsection{One-dimensional scalar problem}\label{SUBSEC:1D_SCALAR}
{\color{orange}We consider a one-dimensional scalar steady state problem with a source term}
\begin{equation}\label{EQ:SSP1D}
f(u)_{x} = s(u, x).
\end{equation}
{\color{orange}On the uniform grid $\left\{x_{i}\right\}_{i=0,\cdots, N}$ with constant $\Delta x = x_{i+1}-x_{i}$, we define the grid function to be $\left\{u_{i}\right\}_{i=0, \cdots, N}$,} the interval $I_{i+\frac{1}{2}}=\left[x_{i}, x_{i+1}\right]$, the control volume $C_{i}=[x_{i-\frac{1}{2}}, x_{i+\frac{1}{2}}]$ and $x_{i+\frac{1}{2}}=\frac{x_{i}+x_{i+1}}{2}$, and the length of $C_{i}$ is denoted by $\left|C_{i}\right|$, which is equal to $\Delta x$.

The total residual in the interval $I_{i+\frac{1}{2}}$ is defined by
\begin{equation}\label{DEF:RES1D}
\Phi_{i+\frac{1}{2}} = \int^{x_{i+1}}_{x_{i}}\!(f(u)_{x} - s(u, x))\,\mathrm{d}x = f(u_{i+1}) - f(u_{i}) - \int^{x_{i+1}}_{x_{i}}\! s(u, x)\, \mathrm{d}x.
\end{equation}
{\color{orange}If we can reach a zero residual limit,} i.e. $\Phi_{i+\frac{1}{2}}=0$ for all $i$, then the accuracy of the scheme is determined by the accuracy of the approximation to $\int^{x_{i+1}}_{x_{i}}\! s(u,x)\, \mathrm{d}x$. In our scheme, we develop a sixth order HWENO integration to approximate the integral $\int^{x_{i+1}}_{x_{i}}\!s(u, x)\,\mathrm{d}x$. {\color{orange}Note that the HWENO integration involves both the point value of the solution $u$  and its spatial derivative $u_{x}$ denoted as $v$, and we denote the corresponding grid function as $\left\{v_{i}\right\}_{i=0, \cdots, N}$.} We now explain the procedure of the sixth order HWENO integration.
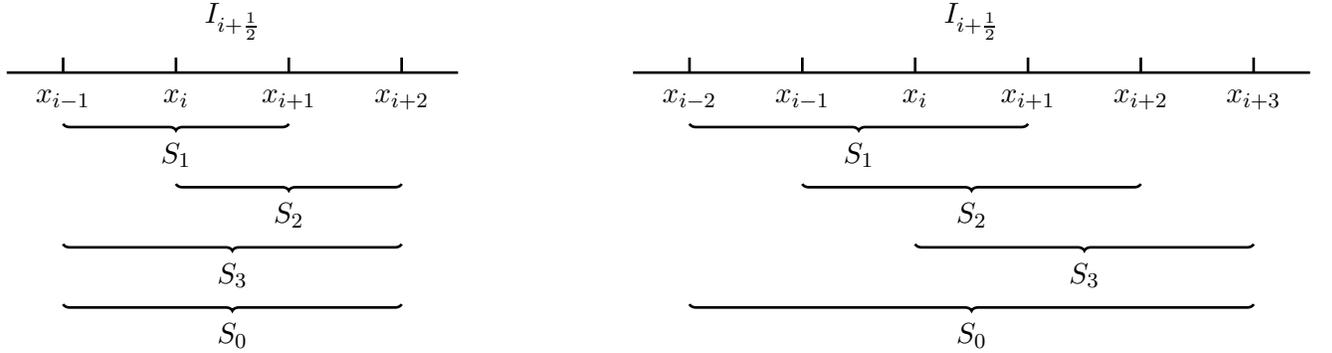
\begin{figure}[!htbp]
\centering
\begin{minipage}{0.45\textwidth}
\begin{tikzpicture}

	\draw[line width=1pt] (0.75, 0)--(6.75, 0);
	
	\draw[xshift = 1.5cm, line width=1pt] (0, 0)--(0, 0.2);
	\node[below] at (1.5, -0.1) {$x_{i-1}$};

	\draw[xshift = 3cm, line width=1pt] (0, 0)--(0, 0.2);
	\node[below] at (3, -0.1) {$x_{i}$};

	\node[above] at (3.75, 0.3) {$I_{i+\frac{1}{2}}$};
	
	\draw[xshift = 4.5cm, line width=1pt] (0, 0)--(0, 0.2);
	\node[below] at (4.5, -0.1) {$x_{i+1}$};

	\draw[xshift = 6cm, line width=1pt] (0, 0)--(0, 0.2);
	\node[below] at (6, -0.1) {$x_{i+2}$};
	
	\draw[decorate,decoration={brace, mirror, raise=8pt}, line width=1pt](1.5, -0.4)--(4.5, -0.4);
	\node[below] at (3, -0.8) {$S_{1}$};
	
	\draw[decorate,decoration={brace, mirror, raise=8pt}, line width=1pt](3, -1.2)--(6, -1.2);
	\node[below] at (4.5, -1.6) {$S_{2}$};
	
	\draw[decorate,decoration={brace, mirror, raise=8pt}, line width=1pt](1.5, -2)--(6, -2);
	\node[below] at (3.75, -2.4) {$S_{3}$};

	\draw[decorate,decoration={brace, mirror, raise=8pt}, line width=1pt](1.5, -2.8)--(6, -2.8);
	\node[below] at (3.75, -3.2) {$S_{0}$};
\end{tikzpicture}
\end{minipage}
\qquad
\begin{minipage}{0.45\textwidth}
\begin{tikzpicture}
	\draw[line width=1pt] (0.75, 0)--(9.75, 0);
	
	\draw[xshift = 1.5cm, line width=1pt] (0, 0)--(0, 0.2);
	\node[below] at (1.5, -0.1) {$x_{i-2}$};
	
	\draw[xshift = 3cm, line width=1pt] (0, 0)--(0, 0.2);
	\node[below] at (3, -0.1) {$x_{i-1}$};

	\draw[xshift = 4.5cm, line width=1pt] (0, 0)--(0, 0.2);
	\node[below] at (4.5, -0.1) {$x_{i}$};
	
	\node[above] at (5.25, 0.3) {$I_{i+\frac{1}{2}}$};
	
	\draw[xshift = 6cm, line width=1pt] (0, 0)--(0, 0.2);
	\node[below] at (6, -0.1) {$x_{i+1}$};

	\draw[xshift = 7.5cm, line width=1pt] (0, 0)--(0, 0.2);
	\node[below] at (7.5, -0.1) {$x_{i+2}$};
	
	\draw[xshift = 9cm, line width=1pt] (0, 0)--(0, 0.2);
	\node[below] at (9, -0.1) {$x_{i+3}$};
	
	\draw[decorate,decoration={brace, mirror, raise=8pt}, line width=1pt](1.5, -0.4)--(6, -0.4);
	\node[below] at (3.75, -0.8) {$S_{1}$};
	
	\draw[decorate,decoration={brace, mirror, raise=8pt}, line width=1pt](3, -1.2)--(7.5, -1.2);
	\node[below] at (5.25, -1.6) {$S_{2}$};
	
	\draw[decorate,decoration={brace, mirror, raise=8pt}, line width=1pt](4.5, -2)--(9, -2);
	\node[below] at (6.75, -2.4) {$S_{3}$};
	
		\draw[decorate,decoration={brace, mirror, raise=8pt}, line width=1pt](1.5, -2.8)--(9, -2.8);
	\node[below] at (5.25, -3.2) {$S_{0}$};
\end{tikzpicture}	
\end{minipage}
%
%
%
%
%
%
%
%
%
\caption{{\color{blue}Stencils used in the sixth order HWENO6 integration (left) vs. WENO6 integration (right).}}
\label{Fig:Stencil_HWENO6-WENO6}
\end{figure}

\begin{description}
\item[Step 1.] Choose a big stencil $S_{0} = \{x_{i-1}, x_{i}, x_{i+1}, x_{i+2}\}$, and construct a fifth degree polynomial $p_{0}(x)$ on it, satisfying
\begin{equation}
\left\{\begin{array}{lll}
p_{0}(x_{i+l}) &=& s_{i+l}, ~~l=-1, 0, 1, 2,\\
p_{0}'(x_{i+l}) &=& s'_{i+l}, ~~l= 0, 1.
\end{array}\right.
\end{equation}
And also choose three small stencils $S_{1}=\left\{x_{i-1}, x_{i}, x_{i+1}\right\}$, $S_{2}=\left\{x_{i}, x_{i+1}, x_{i+2}\right\}$ and $S_{3}=\left\{x_{i-1}, x_{i}, x_{i+1}, x_{i+2}\right\}$, and construct cubic polynomials $p_{1}$, $p_{2}$ and $p_{3}$ on these stencils respectively, which satisfy
\begin{equation}
\left\{\begin{array}{lll}
p_{1}(x_{i+l}) &=& s_{i+l}, ~~l=-1, 0, 1,\\
p'_{1}(x_{i}) &=& s'_{i},
\end{array}\right.
\end{equation}

{\color{orange}
\begin{equation}
\left\{\begin{array}{lll}
p_{2}(x_{i+l}) &=& s_{i+l}, ~~l=0, 1, 2,\\
p'_{2}(x_{i+1}) &=& s'_{i+1},
\end{array}\right.
\end{equation}
}
and
{\color{orange}
\begin{equation}
p_{3}(x_{i+l}) = s_{i+l}, ~~l=-1, 0, 1, 2,
\end{equation}
}
where
\begin{equation*}
\left\{\begin{split}
s_{i+l} &=s(u_{i+l}, x_{i+l}), ~~l=-1, 0, 1, 2,\\
s'_{i+l} &=s'(u_{i+l}, x_{i+l})=\left.\left(\frac{\partial s}{\partial u}v+\frac{\partial s}{\partial x}\right)\right|_{x=x_{i+l}}, ~~l=0, 1.
\end{split}\right.
\end{equation*}
{\color{blue}The chosen stencils between the sixth order HWENO integration (denoted HWENO6) and the sixth order WENO integration (denoted WENO6) are shown in Figure \ref{Fig:Stencil_HWENO6-WENO6}, which indicates that HWENO scheme is more compact than WENO scheme, when designing the same order of accuracy.} And then integrate the above polynomials $p_{0}(x)$, $p_{1}(x)$, $p_{2}$, $p_{3}(x)$ over the interval $I_{i+\frac{1}{2}}$, {\color{orange}and denote them as $q_{0}$, $q_{1}$, $q_{2}$, $q_{3}$, respectively.} In particular, we have
{\color{orange}
\begin{align*}
\begin{split}
q_{0} &=\frac{\Delta x}{240}(s_{i-1}+119s_{i}+119s_{i+1}+s_{i+2}+22\Delta x s'_{i}-22\Delta x s'_{i+1}), \\
q_{1}&= \frac{\Delta x}{24}(s_{i-1}+16s_{i}+7s_{i+1}+6\Delta x s'_{i}),\\
q_{2}&= -\frac{\Delta x}{24}(-7s_{i}-16s_{i+1}-s_{i+2}+6\Delta x s'_{i+1}),\\
q_{3}&= \frac{\Delta x}{24}(-s_{i-1}+13s_{i}+13s_{i+1}-s_{i+2}).\\
\end{split}
\end{align*}
}

\item[Step 2.] Compute the combination coefficients of $q_{l}, l = 1, 2, 3$, such that
\begin{equation}
q_{0} = \sum\limits^{3}_{l=1}\gamma_{l}q_{l},
\end{equation}
{\color{orange}which are so-called linear weights.} Therefore, we obtain
\begin{equation}
\gamma_{1} = \frac{11}{30}, ~~\gamma_{2} = \frac{11}{30}, ~~ \gamma_{3} = \frac{4}{15}.
\end{equation}

\item[Step 3.] Compute the smoothness indicators $\beta_{l}, l=1, 2, 3$, which measure how smooth the functions $p_{l}, l=1, 2, 3$ are in the target interval $I_{i+\frac{1}{2}}$. The smaller these smoothness indicators, the smoother the functions are in the interval $I_{i+\frac{1}{2}}$. {\color{orange}We use the same recipe for the smoothness indicators as in \cite{Balsara.Rumpf_Dumbser.Munz_JCP2009, Jiang.Shu_JCP1996, Shu_SIAMRev2009}:}
\begin{equation}
\beta_{l} = \sum\limits^{r_{0}}_{m=1}\int\limits_{I_{i+\frac{1}{2}}}\Delta x^{2m-1}\!\left(\frac{d^{m}p_{l}}{dx^{m}}\right)^{2}\, \mathrm{d}x, ~~l=1, 2, 3,
\end{equation}
where $r_{0} = 3$ is the corresponding degree of the polynomial. In particular, we have
\begin{align*}
\beta_{1} &= \frac{301}{30}(s_{i-1}+\frac{65}{602}s_{i}-\frac{667}{602}s_{i+1}+\frac{1269}{602}\Delta x
s'_{i})^{2} \\
&~~+\frac{12561}{2408}(s_{i}-s_{i+1}+\frac{10153}{12561}\Delta x s'_{i})^{2}+\frac{10153}{12561}(\Delta x s'_{i})^{2}, \\
\beta_{2} &= \frac{301}{30}(-\frac{667}{602}s_{i}+\frac{65}{602}s_{i+1}+s_{i+2}-\frac{1269}{602}\Delta x s'_{i+1})^{2} \\
&~~+\frac{12561}{2408}(-s_{i}+s_{i+1}-\frac{10153}{12561}\Delta x s'_{i+1})^{2}+\frac{10153}{12561}(\Delta x s'_{i+1})^{2}, \\
\beta_{3} &= \frac{61}{45}(s_{i-1}-\frac{1269}{488}s_{i}+\frac{537}{244}s_{i+1}-\frac{293}{488}s_{i+2})^{2} \\
&~~+\frac{21865}{11712}(s_{i}-\frac{32018}{21865}s_{i+1}+\frac{10153}{21865}s_{i+2})^{2}+\frac{10153}{21865}(s_{i+1}-s_{i+2})^{2}.
\end{align*}

\item[Step 4.] Calculate the non-linear weights based on the linear weights and the smoothness indicators. They are defined as follows:
\begin{equation}
\omega_{l} = \frac{\bar{\omega}_{l}}{\sum^{3}_{m=1}\bar{\omega}_{m}},~~\bar{\omega}_{k} = \frac{\gamma_{l}}{(\varepsilon+\beta_{l})^{2}}, ~~l=1, 2, 3.
\end{equation}

{\color{orange}And it is easy to verify that $\omega_{l} = \gamma_{l}+\mathcal{O}(\Delta x^{3})$ when the solution is smooth.} In our numerical experiments, we take $\varepsilon = 10^{-10}$.

\item[Step 5.] The sixth order reconstruction of the integral $\int^{x_{i+1}}_{x_{i}}\!s(u, x)\, \mathrm{d}x$ is obtained by
\begin{equation}
\int^{x_{i+1}}_{x_{i}}\!s(u, x)\, \mathrm{d}x = \sum\limits^{3}_{m=1}\omega_{m}q_{m}+\mathcal{O}(\Delta x^{7}),
\end{equation}
{\color{orange}when the solution is smooth.}
\end{description}

\begin{rem}
The sixth order HWENO integration leads to the seventh order HWENO approximation to the integral $\int^{x_{i+1}}_{x_{i}}\!s(u, x)\,\mathrm{d}x$ within each interval and hence the sixth order approximation to the integral over the whole computational domain.
\end{rem}

Next we start to distribute the total residuals. In the interval $[x_{i}, x_{i+1}]$,  the total residual is $\Phi_{i+\frac{1}{2}}$, {\color{orange}and it is distributed to the nodes $x_{i}$ and $x_{i+1}$ denoted as $\Phi_{i+\frac{1}{2}}^{1}$ and $\Phi_{i+\frac{1}{2}}^{2}$.} For simplicity and with no ambiguity, we drop off the subscript $i+\frac{1}{2}$ for the total residual $\Phi_{i+\frac{1}{2}}$. We require that $\Phi=\Phi^{1}+\Phi^{2}$ for the conservation and $|\Phi^{k}|/|\Phi|, k=1, 2$ be uniformly bounded by the residual property \cite{Abgrall_JCP2001}. {\color{orange}One way to distribute the total residual $\Phi$ with the upwinding property is given by} 
\begin{equation}\label{UW_RD_HWENO1D}
\Phi^{1} = (1-\alpha)\Phi, ~~\Phi^{2} = \alpha\Phi, ~~\alpha\in [0, 1],
\end{equation}
with $\alpha$ defined by
\begin{equation*}
\alpha = \left\{\begin{array}{ll}
1 & \mbox{if}~\bar{\lambda} \geq \delta, \\
0 & \mbox{if}~\bar{\lambda} \leq -\delta, \\
r(\bar{\lambda}, \delta) & \mbox{otherwise},
\end{array} \right.
\end{equation*}
where $\bar{\lambda} = f'(\bar{u})$, {\color{orange}and $\bar{u}$ is an average state in the interval taken to be $\frac{1}{2}(u_{i}+u_{i+1})$.} The function $r(\cdot, \cdot)$ is a continuous differentiable entropy function for the Roe scheme \cite{Harten_JCP1983}, which is given by
\begin{equation}\label{ROE_FUN}
r(\lambda, \delta) = \frac{1}{4\delta^{3}}(\lambda+\delta)^{2}(2\delta-\lambda),
\end{equation}
where the coefficient $\delta$ is chosen accordingly in the problem.

{\color{orange}Let $u_{i}^{\text{new}}$ and $v_{i}^{\text{new}}$ denote the updated numerical approximation of $u$ and the corresponding derivative $v$, respectively, at the grid point $x_{i}$.  And $u_{i}^{\text{new}}$ is updated through sending the distributed residuals to the point $x_{i}$,} as in a pseudo time-marching scheme, which can be written as a semi-discrete system
\begin{equation}\label{SCHEMA:UW_HWENO_SCALAR1D}
\frac{\mathrm{d}u^{\text{new}}_{i}}{\mathrm{d}t} + \frac{1}{\left|C_{i}\right|}\left(\Phi^{2}_{i-\frac{1}{2}}+\Phi^{1}_{i+\frac{1}{2}}\right)=0.
\end{equation}
In our numerical experiments, we use a forward Euler scheme for the pseudo time discretization. {\color{orange}As mentioned above, a zero residual limit $\Phi_{i+\frac{1}{2}}$ is a steady state solution of \eqref{SCHEMA:UW_HWENO_SCALAR1D}. However, near shocks, $\Phi_{i+\frac{1}{2}}$ may not be very small even if the steady state solution is reached according to our numerical experiments. Here we refer readers to \cite{Chou.Shu_JCP2006} for the convergence towards weak solutions in this situation via a Lax-Wendroff type theorem.}

{\color{orange}
In the framework of the traditional HWENO scheme, we need to take partial derivative w.r.t. the variable $x$ on the both side of the equation \eqref{EQ:SSP1D}, then we get an additional equation related to the spatial derivative $v(x, t)$:
\begin{equation}
f_{1}(u, v)_{x}=s(u, x)_{x},
\end{equation}
where $f_{1}(u, v)=f'(u)u_{x}=f'(u)v$. Following the recipe of the RD scheme for the variable $u(x, t)$, we can also define the total residual w.r.t. $v(x, t)$ through the integral form and distribute the residuals in an upwinding way as the variable $u(x, t)$. Finally, update the derivative of the point value $v^{\text{new}}_{i}$ through sending the distributed residuals to the point $x_{i}$ as in a pseudo time-marching scheme, see \cite{mythesis} for more details. However, in this paper, we use the current point value $u^{\text{new}}_{i}$ and the old spatial derivative $v^{\text{old}}_{i}$ to reconstruct  the spatial derivative of the point value $v^{\text{new}}_{i}$, {\color{orange}due to the fact that this is a steady state problem and it does not involve time, so that $v_{i}^{\text{old}}$ is also a ``good'' approximation to the exact solution}. Hence, we do not need to introduce an additional auxiliary equation, which saves computational storage. Here we adopt a fourth order HWENO reconstruction proposed in \cite{ZHao.Zhang.Qiu_JSC2020}. We now roughly recall the procedure of the reconstruction. For simplicity, {\color{orange}in the following reconstruction process, we drop off the superscript ``new'' for $u^{\text{new}}$ and the superscript ``old'' for $v^{\text{old}}$}.}

\begin{figure}[!htbp]
\centering
\begin{tikzpicture}
\draw[line width=1pt] (0.75, 0)--(5.25, 0);
\draw[xshift = 1.5cm, line width=1pt] (0, 0)--(0, 0.2);
\node[below] at (1.5, -0.1) {$x_{i-1}$};
\node[above] at (1.5, 0.3) {$u_{i-1}$};
\node[above] at (1.5, 0.8) {$v_{i-1}$};

\draw[xshift = 3cm, line width=1pt] (0, 0)--(0, 0.2);
\node[below] at (3, -0.1) {$x_{i}$};
\node[above] at (3, 0.3) {$u_{i}$};
\node[fill, circle, scale = 0.5, blue] at (3, 0) {};

\draw[xshift = 4.5cm, line width=1pt] (0, 0)--(0, 0.2);
\node[below] at (4.5, -0.1) {$x_{i+1}$};
\node[above] at (4.5, 0.3) {$u_{i+1}$};
\node[above] at (4.5, 0.8) {$v_{i+1}$};

\draw[decorate,decoration={brace, mirror, raise=8pt}, line width=1pt](1.5, -0.4)--(3, -0.4);
\node[below] at (2.25, -0.8) {$S_{1}$};
\draw[decorate,decoration={brace, mirror, raise=8pt}, line width=1pt](3, -1.2)--(4.5, -1.2);
\node[below] at (3.75, -1.6) {$S_{2}$};
\draw[decorate,decoration={brace, mirror, raise=8pt}, line width=1pt](1.5, -2)--(4.5, -2);
\node[below] at (3, -2.4) {$S_{3}$};

\draw[decorate,decoration={brace, mirror, raise=8pt}, line width=1pt](1.5, -2.8)--(4.5, -2.8);
\node[below] at (3, -3.2) {$S_{0}$};
\end{tikzpicture}
\caption{{\color{orange}Stencils used in the fourth order HWENO reconstruction.}}
\label{PIC:stencil_hweno4}
\end{figure}
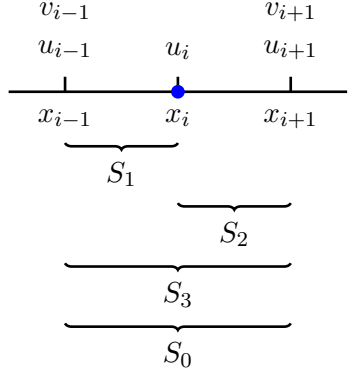

\begin{description}
\item[{\color{orange}Step 1.}] {\color{orange}Choose three stencils $S_{1} = \{x_{i-1}, x_{i}\}$, $S_{2} = \{x_{i}, x_{i+1}\}$, $S_{3} = \{x_{i-1}, x_{i}, x_{i+1}\}$, and a big stencil $S_{0}=\left\{x_{i-1}, x_{i}, x_{i+1}\right\} $, as shown in Figure \ref{PIC:stencil_hweno4}. Then construct three quadratic polynomials $P_{1}(x)$, $P_{2}(x)$ and $P_{3}(x)$ on the stencils $S_{1}$, $S_{2}$ and $S_{3}$ by the Hermite interpolation,} respectively, under the following conditions
\begin{align*}
P_{1}(x_{i+l}) &= u_{i+l}, ~l = -1, 0, ~ P'_{1}(x_{i-1}) = v_{i-1},\\
P_{2}(x_{i+l}) &= u_{i+l}, ~l = 0, 1, ~P'_{2}(x_{i+1}) = v_{i+1},\\
P_{3}(x_{i+l}) &= u_{i+l}, ~l=-1, 0, 1.
\end{align*}
Similarly, {\color{orange}a quartic polynomial $P_{0}$ is also obtained by the Hermite interpolation,} satisfying
\begin{equation*}
P_{0}(x_{i+l}) = u_{i+l}, ~l=-1, 0, 1, ~P_{0}'(x_{i+l})=v_{i+l}, ~l = -1, 1.
\end{equation*}
Therefore, {\color{orange}we get their derivative values at the point $x_{i}$, respectively in the following:}
\begin{align*}
P'_{1}(x_{i}) &= -v_{i-1} + \frac{2(u_{i}-u_{i-1})}{\Delta x}, \\
P'_{2}(x_{i}) &= -v_{i+1}+\frac{2(u_{i+1}-u_{i})}{\Delta x}, \\
P'_{3}(x_{i}) &= \frac{u_{i+1}-u_{i-1}}{\Delta x}, \\
P_{0}'(x_{i}) &= -\frac{v_{i-1}+v_{i+1}}{4}+\frac{3(u_{i+1}-u_{i-1})}{4\Delta x}.
\end{align*}
\item[Step 2.] Compute the linear weights $\gamma_{l}, l = 1, 2, 3$ by requiring that
\begin{equation*}
P_{0}'(x_{i}) = \sum\limits^{3}_{l=1}\gamma_{l}P'_{l}(x_{i}),
\end{equation*}
then we have
\begin{align*}
\gamma_{1} = \frac{1}{4},~\gamma_{2} = \frac{1}{4}, ~\gamma_{3} = \frac{1}{2}.
\end{align*}

\item[Step 3.] Compute the smoothness indicators $\beta_{l}, l = 1, 2, 3$ to measure how smooth the functions $P_{l}(x), l = 1, 2, 3$ are in the target cell $C_{i}$, defined in the following
{\color{orange}
\begin{equation}
\beta_{l} = \sum\limits^{r_{0}}_{m=1}\int\limits_{C_{i}}\Delta x^{2m-1}\!\left(\frac{d^{m}P_{l}}{dx^{m}}\right)^{2}\, \mathrm{d}x,
\end{equation}
where $r_{0} = 2$ is the degree of the polynomials $P_{l}, l=1, 2, 3$.} In particular, we have
\begin{align*}
\beta_{1} &= (2u_{i}-2u_{i-1}-v_{i-1}\Delta x)^{2}+\frac{13}{3}(u_{i}-u_{i-1}-v_{i-1}\Delta x)^{2}, \\
\beta_{2} &= (2u_{i+1}-2u_{i}-v_{i+1}\Delta x)^{2}+\frac{13}{3}(u_{i+1}-u_{i}-v_{i+1}\Delta x)^{2}, \\
\beta_{3} &= \frac{1}{4}(u_{i+1}-u_{i-1})^{2}+\frac{13}{12}(u_{i-1}-2u_{i}+u_{i+1})^{2}.
\end{align*}

\item[Step 4.] {\color{orange}The non-linear weights are obtained as}
\begin{equation}
\omega_{l} = \frac{\bar{\omega}_{l}}{\sum^{3}_{m=1}\bar{\omega}_{m}}, ~~\bar{\omega}_{l} = \frac{\gamma_{l}}{(\beta_{l}+\varepsilon)^{2}},~~l = 1, 2, 3,
\end{equation}
where $\varepsilon$ is taken to be $10^{-10}$ to avoid denominator being zero. Therefore, the derivative value $v^{\text{new}}_{i}$ is obtained by
{\color{orange}
\begin{equation}
v^{\text{new}}_{i} = \sum\limits^{3}_{l=1}\omega_{l}P'_{l}(x_{i})+\mathcal{O}(\Delta x^{4}).
\end{equation}
}
\end{description}

{\color{blue}
\begin{rem}
From the descriptions of the novel HWENO framework and traditional HWENO framework above, it is the same process to deal with the unknown variable $u$ for both schemes. Hence, the novel HWENO scheme keeps the same compactness of stencil as the traditional HWENO framework.
\end{rem}	
}

{\color{blue}
\begin{rem}
For the sake of the compactness of stencils and good performances of the scheme, we employ the fourth order HWENO reconstruction. Actually, the expected accuracy of the spatial derivative $v^{\text{new}}_{i}$ is the fifth order accuracy, which can be complied with the sixth order HWENO integration. In our error and order of accuracy of numerical experiments indicate that the RD conservative finite difference HWENO scheme is the sixth order accuracy at least in the $L^{1}$ sense, when adopting the fourth order HWENO reconstruction for the derivative $v^{\text{new}}_{i}$.
\end{rem}
}

{\color{orange}
\underline{\bf The procedure \rom{1}. RD finite difference HWENO method for steady state problem in 1D}
\begin{enumerate}
\item Compute the total residual defined in \eqref{DEF:RES1D}, using the sixth order HWENO integration to approximate the source term.

\item Distribute the total residual within each interval to two nodes $x_{i}$ and $x_{i+1}$, complying with upwinding principle given in  \eqref{UW_RD_HWENO1D}.

\item Update the point value $u^{\text{new}}_{i}$ through sending the distributed residuals to the point $x_{i}$ and forward in a pseudo time time-marching scheme by a forward Euler time discretization until the steady state is reached.

\item Update the spatial derivative of the point value $v^{\text{new}}_{i}$, followed by the fourth HWENO reconstruction, using the current point value $u^{\text{new}}_{i}$ and the old spatial derivative of the point value $v^{\text{old}}_{i}$ to reconstruct.
\end{enumerate}
}

\subsection{One-dimensional systems}
Consider a one-dimensional steady state system \eqref{EQ:SSP1D}, where {\bf u}, {\bf f}({\bf u}) and {\bf s}({\bf u}, x) are vector-valued functions in $\mathbb{R}^{m}$. For hyperbolic systems, we assume that the Jacobian matrix ${\bf f}'({\bf u})$ can be written as $R\Lambda L$, where $\Lambda$ is a diagonal matrix with real eigenvalues on the diagonal, {\color{orange}and $L$ and $R$ are matrices of left and right eigenvectors of ${\bf f}'({\bf u})$, respectively.}

The grid, grid function, {\color{orange}grid function of the spatial derivative,} the interval, and control volume are denoted as in Subsection \ref{SUBSEC:1D_SCALAR}. The total residual ${\bf \Phi}_{i+\frac{1}{2}} $ in the interval $\left[x_{i}, x_{i+1}\right]$ is again defined by \eqref{DEF:RES1D}. As before, the accuracy of the scheme is determined by the accuracy of the approximation to $\int^{x_{i+1}}_{x_{i}}\!{\bf s}({\bf u}, x)\, \mathrm{d}x$, {\color{orange}which is again approximated by the sixth order HWENO integration introduced in Subsection \ref{SUBSEC:1D_SCALAR}.}

In order to distribute the total residual ${\bf \Phi}_{i+\frac{1}{2}}$, we need to use a local characteristic decomposition in the interval $[x_{i}, x_{i+1}]$. First, we compute an average state $\bar{{\bf u}}$ between ${\bf u}_{i}$ and ${\bf u}_{i+1}$, using either the simple arithmetic mean or Roe's average \cite{Roe_JCP1981}, {\color{orange} and $\bar{L}$ and $\bar{R}$ are the corresponding left and right eigenvectors $L$ and $R$ evaluated at the average state $\bar{{\bf u}}$,} and $\bar{\lambda}_{k}$ is the corresponding $k$-th eigenvalue. In the following, for simplicity and with no ambiguity, we drop off the subscript $i+\frac{1}{2}$ for the total residual ${\bf \Phi}_{i+\frac{1}{2}}$. To keep the conservation, we require that ${\pmb \Phi}={\pmb \Phi}^{1}+{\pmb \Phi}^{2}$ and $|{\pmb \Phi}^{k}|/|{\pmb \Phi}|, k=1, 2$, which is in a component by component sense and similar to the one-dimensional scalar cases, be uniformly bounded to guarantee the residual property \cite{Abgrall_JCP2001}. Now we project total residual ${\pmb \Phi}$ at the interval to the characteristic field, namely, ${\pmb \Psi}=\bar{L}{\pmb \Phi}$, then distribute the residual ${\pmb \Psi}$ to the nodes $x_{i}$ and $x_{i+1}$ in the upwinding way, denoted by ${\pmb \Psi}^{1}$ and ${\pmb \Psi}^{2}$, and we obtain that ${\pmb \Psi}={\pmb \Psi}^{1}+{\pmb \Psi}^{2}$. According to one-dimensional scalar problem, the upwind scheme can be defined naturally in the following:
\begin{equation}\label{UW_RD_HWENO1D_SYS}
{\pmb \Psi}^{1}=(\mathbf{I}-\Sigma){\pmb \Psi}, ~~{\pmb \Psi}^{2}=\Sigma{\pmb \Psi},
\end{equation}
where $\mathbf{I}$ is an identity matrix, $\Sigma$ is a diagonal matrix, and the $k$-th diagonal component of $\Sigma$ is given by
\begin{equation*}
\Sigma_{kk}=\left\{\begin{array}{ll}
1 & \mbox{if}~\bar{\lambda}_{k} \geq \delta, \\
0 & \mbox{if}~\bar{\lambda}_{k} \leq -\delta, \\
r(\bar{\lambda}_{k}, \delta) & \mbox{otherwise}.
\end{array}\right.
\end{equation*}
{\color{orange}Here the function $r(\cdot, \cdot)$ is defined as \eqref{ROE_FUN},} and $\delta$ is also chosen accordingly in the problem.

{\color{orange}Next we need to project the residuals ${\pmb \Psi}^{1}$ and ${\pmb \Psi}^{2}$ back to the physical space and then obtain the residuals distributed to the nodes $x_{i}$ and $x_{i+1}$ as follows:}
\begin{equation}\label{UW_RD_HWENO1D_SYS_pbk}
{\pmb \Phi}^{1}=\bar{R}{\pmb \Psi}^{1}, ~~{\pmb \Phi}^{2}=\bar{R}{\pmb \Psi}^{2}.
\end{equation}

Thus, we get the way to distribute the total residual within each interval. As in the scalar case, the point value ${\bf u}^{\text{new}}_{i}$ can be updated in the pseudo time-marching semi-discrete scheme \eqref{SCHEMA:UW_HWENO_SCALAR1D}, which is again discretized by a forward Euler scheme in our numerical experiments until the steady state is reached. {\color{orange}As for ${\bf v}^{\text{new}}_{i}$, it is again reconstructed component by component via the fourth order HWENO reconstruction as described in Subsection \ref{SUBSEC:1D_SCALAR}.}

{\color{orange}
\underline{\bf Procedure \rom{2}. RD finite difference HWENO method for steady state problem in 1D}
\begin{enumerate}
\item Compute the total residual as defined in \eqref{DEF:RES1D} component by component, using the sixth order HWENO integration to approximate the source term.

\item Project the total residual to a local characteristic field, and then distribute it to the nodes $x_{i}$ and $x_{i+1}$  in the upwinding way, as given in \eqref{UW_RD_HWENO1D_SYS}; then project the residuals back to physical space, as in \eqref{UW_RD_HWENO1D_SYS_pbk}.

\item Update the point value ${\bf u}^{\text{new}}_{i}$ through sending the distributed residuals to the point $x_{i}$ in the physical space and forward in a pseudo time \eqref{SCHEMA:UW_HWENO_SCALAR1D} by a forward Euler time discretization until the steady state is reached.

\item Update the spatial derivative of the point value ${\bf v}^{\text{new}}_{i}$, which is reconstructed component by component by the fourth order HWENO scheme, using the current point value ${\bf u}^{\text{new}}_{i}$ and the old spatial derivative of the point value ${\bf v}^{\text{old}}_{i}$ to reconstruct.
\end{enumerate}
}

\section{High order RD conservative finite difference HWENO scheme in two dimensions}
In this section, {\color{orange}we develop a high order RD conservative finite difference HWENO scheme for two-dimensional steady state problems.} More precisely, {\color{orange} we focus on our scheme on uniform Cartesian meshes.} In Subsection \ref{SUBSEC:2D_SCALAR}, we define the total residual within each cell through the integral form, as defined in \eqref{DEF:RES1D}, and then introduce the residual distribution mechanism. In Subsection \ref{SUBSEC:2D_SYS}, {\color{red} we extend the scheme to two-dimensional systems, which is based on a local characteristic field decomposition, and  distribute the total residuals in characteristic fields dimension-by-dimension}.

\subsection{Two-dimensional scalar problems}\label{SUBSEC:2D_SCALAR}

{\color{orange}We consider a two-dimensional scalar steady state problem with a source term}
\begin{equation}\label{EQ:SSP2D}
f(u)_{x}+g(u)_{y}=s(u, x, y).
\end{equation}
{\color{orange}On the uniform grid to be $\left\{(x_{i}, y_{j})\right\}_{i=0, \cdots, N \atop j=0, \cdots, M}$ with constant $\Delta x = x_{i+1}-x_{i}$ and $\Delta y = y_{j+1}-y_{j}$, we define the grid function to be $u_{i,j}$ ,} the cell $I_{i+\frac{1}{2}, j+\frac{1}{2}}=\left[x_{i}, x_{i+1}\right]\times \left[y_{j}, y_{j+1}\right]$, the control volume $C_{ij}=[x_{i-\frac{1}{2}}, x_{i+\frac{1}{2}}]\times[y_{j-\frac{1}{2}}$, $y_{j+\frac{1}{2}}]$, and $x_{i+\frac{1}{2}}=\frac{x_{i}+x_{i+1}}{2}$ and $y_{j+\frac{1}{2}}=\frac{y_{j}+y_{j+1}}{2}$, and the area of $C_{ij}$ is denoted by $\left|C_{ij}\right|$, which is equal to $\Delta x \cdot \Delta y$.

The total residual in the cell $I_{i+\frac{1}{2}, j+\frac{1}{2}}$ is defined by
\begin{align}\label{DEF:RES2D}
\begin{split}
\Phi_{i+\frac{1}{2},j+\frac{1}{2}} & =  \int^{y_{j+1}}_{y_{j}}\int^{x_{i+1}}_{x_{i}}\!(f(u)_{x} + g(u)_{y} - s(u, x, y))\,\mathrm{d}x\,\mathrm{d}y  \\
& = \int^{y_{j+1}}_{y_{j}}\!(f(u(x_{i+1}, y)) - f(u(x_{i}, y)))\,\mathrm{d}y \\
& ~~ + \int^{x_{i+1}}_{x_{i}}\!(g(u(x, y_{j+1})) - g(u(x, y_{j})))\,\mathrm{d}x  \\
& ~~ -\int^{y_{j+1}}_{y_{j}}\int^{x_{i+1}}_{x_{i}}\!s(u(x, y), x, y)\,\mathrm{d}x\,\mathrm{d}y.
\end{split}
\end{align}
{\color{orange}If we can reach a zero residual limit, i.e., $\Phi_{i+\frac{1}{2}, j+\frac{1}{2}} = 0$ for all $i$ and $j$,} the accuracy of the scheme is determined by the accuracy of the approximations to the integrations of the fluxes and the source term.

To approximate the integrations of the fluxes, which are one-dimensional integrals, {\color{orange}we just use the sixth order HWENO integration as described in Subsection \ref{SUBSEC:1D_SCALAR}.} As for the source term $\int^{y_{j+1}}_{y_{j}}\int^{x_{i+1}}_{x_{i}}\!s(u, x, y)\, \mathrm{d}x\, \mathrm{d}y$, {\color{orange} we can approximate it in a dimension-by-dimension fashion, due to the finite difference scheme on the Cartesian meshes,} which is explained as follows.

First, we define
{\color{orange}
\begin{align*}
S_{i+1/2}(y) &= \int^{x_{i+1}}_{x_{i}}\!s(u(x, y), x, y)\,\mathrm{d}x, \\
\left(S_{i+1/2}(y)\right)_{y} &= \int^{x_{i+1}}_{x_{i}}\!  s(u(x, y), x, y)_{y}\,\mathrm{d}x,
\end{align*}
}
and then the integral of the source term in the cell $I_{i+\frac{1}{2}, j+\frac{1}{2}}$ can be rewritten as

{\color{orange}
\begin{equation*}
\int^{y_{j+1}}_{y_{j}}\int^{x_{i+1}}_{x_{i}}\!s(u(x, y), x, y)\,\mathrm{d}x\,\mathrm{d}y = \int^{y_{j+1}}_{y_{j}}\!S_{i+1/2}(y)\,\mathrm{d}y.
\end{equation*}	
The integral $\int^{y_{j+1}}_{y_{j}}\!S_{i+1/2}(y)\,\mathrm{d}x$ can be approximated by the sixth order HWENO integration in the $y$-direction, using $\left\{S_{i+1/2}(y_{j+k})\right\}_{k=-1,\cdots,2}$ and $\left\{\left( S_{i+1/2}(y_{j+k})\right)_{y} \right\}_{k = 0, 1}$. By the definition of $S_{i+1/2}(y)$ and $\left(S_{i+1/2}(y)\right)_{y}$, $S_{i+1/2}(y_{j+k})$ can again be approximated by the sixth order HWENO integration in the $x$-direction,} using $\left\{s(u_{i+l,j+k}, x_{i+l}, y_{j+k})\right\}_{l=-1,\cdots,2}$ and $\left\{ s(u_{i+l, j+k}, x_{i+l}, y_{j+k})_{x}\right\}_{l = 0, 1}$. Similarly,  {\color{orange}$\left(S_{i+1/2}(y_{j+k})\right)_{y}$ can be approximated by the sixth order HWENO integration in the $x$-direction,} using $\left\{s(u_{i+l, j+k}, x_{i+l}, y_{j+k})_{y}\right\}_{l=-1, \cdots, 2}$ and $\left\{s(u_{i+l, j+k}, x_{i+l}, y_{j+k})_{xy}\right\}_{l=0, 1}$. Thus, {\color{orange}the integration of the source term can be approximated dimension-by-dimension,} and the sixth order accuracy is obtained at the zero residual limit.

\begin{rem}
{\color{orange}According to the approximation for the source term above, the spatial derivatives $(u_{x})_{i, j}$, $(u_{y})_{i, j}$ and the second cross derivative $(u_{xy})_{i, j}$ are still involved in the procedure of HWENO integration, which are denoted by $v_{i, j}$, $w_{i, j}$ and $z_{i, j}$, respectively.}
\end{rem}

Next, we start to distribute the total residuals. In the cell $I_{i+\frac{1}{2}, j+\frac{1}{2}}=\left[x_{i},x_{i+1}\right]\times\left[y_{j},y_{j+1}\right]$, the total residual is $\Phi_{i+\frac{1}{2}, j+\frac{1}{2}}$, and it is to be distributed to the vertices of the cell, which are denoted to be $M_{1} = (x_{i}, y_{j})$, $M_{2} = (x_{i+1}, y_{j})$, $M_{3} = (x_{i}, y_{j+1})$ and $M_{4} = (x_{i+1}, y_{j+1})$. Here we denote the residuals distributed to the vertices $M_{k}$ as $\Phi^{k}_{i+\frac{1}{2}, j+\frac{1}{2}}$, $k=1,2,3,4$. For simplicity and without ambiguity, we drop off the subscript $(i+\frac{1}{2}, j+\frac{1}{2})$ in the notations. For the conservation and the residual property in \cite{Abgrall_JCP2001}, we require that $\Phi = \sum^{4}_{k=1}\Phi^{k}$ and $|\Phi^{k}|/|\Phi|, k=1, \cdots, 4$ be uniformly bounded.

{\color{orange}In order to have the upwinding property, we have an upwind scheme in the following}
\begin{equation}\label{UW_RD_HWENO2D}
\tilde{\Phi}^{1} =(1-\alpha)((1-\beta))\Phi,~\tilde{\Phi}^{2} =\alpha(1-\beta)\Phi,~\tilde{\Phi}^{3}=(1-\alpha)\beta\Phi,~\tilde{\Phi}^{4} =\alpha \beta \Phi,
\end{equation}
with $\alpha, \beta \in [0, 1]$. $\alpha$ is the coefficient for upwinding in the $x$-direction, which is given by
\begin{equation*}
\alpha=\left\{\begin{array}{ll}
1 & \mbox{if}~\bar{\lambda}_{x} \geq \delta, \\
0 & \mbox{if}~\bar{\lambda}_{x} \leq -\delta, \\
r(\bar{\lambda}_{x}, \delta) & \mbox{otherwise},
\end{array}\right.
\end{equation*}
where $\bar{\lambda}_{x}=f'(\bar{u})$, and $\bar{u}$ is an average state in the cell $I_{i+\frac{1}{2}, j+\frac{1}{2}}$ to taken by
$$\bar{u}=\frac{1}{4}(u_{i, j}+u_{i+1, j}+u_{i, j+1}+u_{i+1, j+1}).$$
Similarly, $\beta$ is the coefficient for upwinding in the $y$-direction, which is given by
\begin{equation*}
\beta=\left\{\begin{array}{ll}
1 & \mbox{if}~\bar{\lambda}_{y} \geq \delta, \\
0 & \mbox{if}~\bar{\lambda}_{y} \leq -\delta, \\
r(\bar{\lambda}_{y}, \delta) & \mbox{otherwise},
\end{array}\right.
\end{equation*}
where $\bar{\lambda}_{y}=g'(\bar{u})$. And the function $r(\cdot, \cdot)$ is given as in \eqref{ROE_FUN}, $\delta$ is chosen accordingly in the problem.

Our numerical experiments show that for shock problems, it is necessary to add dissipation term for the stability of the scheme proposed. Thus, we introduce an additional dissipation residual $\Phi^{k}_{\mbox{diss}}$ to the residual $\tilde{\Phi}^{k}$ for each vertex. Here is the definition of dissipation residual in the following:
\begin{align}\label{DEF:UW_HWENO2D_RESIDUAL_DISS}
\begin{split}
\Phi^{1}_{\mbox{diss}} &=\frac{\sigma}{2}\Delta^{3}\left(\frac{u_{i, j}-u_{i+1, j}}{\Delta x}+\frac{u_{i, j}-u_{i, j+1}}{\Delta y}\right), \\
\Phi^{2}_{\mbox{diss}} &=\frac{\sigma}{2}\Delta^{3}\left(\frac{u_{i+1, j}-u_{i, j}}{\Delta x}+\frac{u_{i+1, j}-u_{i+1, j+1}}{\Delta y}\right), \\
\Phi^{3}_{\mbox{diss}} &=\frac{\sigma}{2}\Delta^{3}\left(\frac{u_{i, j+1}-u_{i+1, j+1}}{\Delta x}+\frac{u_{i, j+1}-u_{i, j}}{\Delta y}\right),  \\
\Phi^{4}_{\mbox{diss}} &=\frac{\sigma}{2}\Delta^{3}\left(\frac{u_{i+1, j+1}-u_{i, j+1}}{\Delta x}+\frac{u_{i+1, j+1}-u_{i+1, j}}{\Delta y}\right),
\end{split}
\end{align}
where $\Delta=\max\left(\Delta x, \Delta y\right)$, and $\sigma$ is chosen accordingly in the problem. This dissipation mechanism works well for our numerical experiments, but it may not be the optimal approach, since it has a adjustable coefficient $\sigma$, whose choice for optimal performance seems to be problem dependent.

Thus, we get the way to distribute the total residual within each cell and obtained by
\begin{equation}\label{UW_HWENO2D_RES_FINAL}
\Phi^{k}=\tilde{\Phi}^{k}+\Phi^{k}_{\mbox{diss}}, ~~k=1, \cdots, 4.
\end{equation}
{\color{orange}The point value $u^{\text{new}}_{i,j}$ is then updated through sending the distributed residuals to the point $(x_{i}, y_{j})$,} as in a pseudo time-marching scheme, which can be written as a semi-discrete system
{\color{orange}
\begin{equation}\label{SCHEMA:UW_HWENO2D_SCALAR}
\frac{\mathrm{d}u^{\text{new}}_{i,j}}{\mathrm{d}t}+\frac{1}{|C_{ij}|}\left(\Phi^{1}_{i+\frac{1}{2}, j+\frac{1}{2}}+\Phi^{2}_{i-\frac{1}{2}, j+\frac{1}{2}}+\Phi^{3}_{i+\frac{1}{2}, j-\frac{1}{2}}+\Phi^{4}_{i-\frac{1}{2}, j-\frac{1}{2}}\right)=0.
\end{equation}
}
We again use a forward Euler scheme for the pseudo time discretization. And we may lose the strict residual property after adding dissipation residuals, but note that conservation is still preserved after adding the dissipation since $\sum^{4}_{k=1}\Phi^{k}_{\mbox{diss}}=0$. 

{\color{orange}As for the spatial derivatives $v^{\text{new}}_{i, j}$ and $w^{\text{new}}_{i, j}$, they can be approximated by the fourth order HWENO reconstruction in the $x$-direction and in the $y$-direction, respectively, as introduced in Subsection \ref{SUBSEC:1D_SCALAR}. As for the second cross derivative $z^{\text{new}}_{i, j}$, it is updated dimension-by-dimension.} {\color{orange}We first consider $(u_{x})_{i, j}$, which is approximated by the fourth order HWENO reconstruction in the $x$-direction, and then $((u_{x})_{y})_{i, j}$ again by the fourth order HWENO reconstruction in the $y$-direction. Thus, we get the way to update $z^{\text{new}}_{i, j}$.}

{\color{orange}
\underline{\bf Procedure \rom{3}. RD finite difference HWENO method for steady state problem in 2D}
\begin{enumerate}
\item Compute the total residuals within each cell as defined in \eqref{DEF:RES2D}, 
using the sixth order HWENO integration dimension-by-dimension for the source term.

\item Distribute the total residual to the four vertices of the cell via an upwind scheme given in \eqref{UW_RD_HWENO2D}.

\item Revise the residuals distributed to the four vertices by adding additional dissipation residuals, defined in \eqref{DEF:UW_HWENO2D_RESIDUAL_DISS}.

\item Update the point value $u^{\text{new}}_{i, j}$ through sending the distributed residuals to the point $(x_{i}, y_{j})$ and forward in a pseudo time \eqref{SCHEMA:UW_HWENO2D_SCALAR} by a forward Euler time discretization until the steady state is reached.

\item Update the spatial derivatives of the point value $v^{\text{new}}_{i, j}$ and $w^{\text{new}}_{i, j}$ by the fourth order HWENO reconstruction in the $x$-direction and in the $y$-direction, respectively, as in Subsection \ref{SUBSEC:1D_SCALAR}. The second cross derivative $z^{\text{new}}_{i, j}$ is reconstructed by the fourth order HWENO reconstruction in a dimension-by-dimension way.
\end{enumerate}
}

\subsection{Two-dimensional systems}\label{SUBSEC:2D_SYS}
Consider a two-dimensional steady state system \eqref{EQ:SSP2D}, where {\bf u}, {\bf f}({\bf u}), {\bf g}({\bf u}) and {\bf s}({\bf u}, x, y) are vector-valued functions in $\mathbb{R}^{m}$. For hyperbolic systems, we assume that any real linear combination of the Jacobians $n_{x}{\bf f}'({\bf u}) + n_{y}{\bf g}'({\bf u})$ is diagonalizable with real eigenvalues. In particular, we assume ${\bf f}'({\bf u})$ and ${\bf g}'({\bf u})$ can be written as $R_{x}\Lambda_{x}L_{x}$ and $R_{y}\Lambda_{y}L_{y}$, respectively, where $\Lambda_{x}$ and $\Lambda_{y}$ are diagonal matrices with real eigenvalues on the diagonal, and $L_{x}$, $R_{x}$ and $L_{y}$, $R_{y}$ are matrices of left and right eigenvectors for the corresponding Jacobians.

{\color{orange}The grid, grid function, grid function of the derivatives, cell and control volume are denoted as in Subsection \ref{SUBSEC:2D_SCALAR}.} The total residual in the cell $I_{i+\frac{1}{2}, j+\frac{1}{2}} = \left[x_{i}, x_{i+1}\right]\times\left[y_{j}, y_{j+1}\right]$ is still defined by \eqref{DEF:RES2D}. As before, {\color{orange}if we can reach a zero residual limit of the scheme,} the accuracy of the scheme is determined by the accuracy of the approximations to the integrals of the fluxes and the source term.
{\color{orange}The integrals of the fluxes are again calculated by the sixth order HWENO integration described in Subsection \ref{SUBSEC:1D_SCALAR}, and the integral of the source term is calculated in a dimension-by-dimension fashion by the sixth order HWENO integration as shown in Subsection \ref{SUBSEC:2D_SCALAR}.} For simplicity and without ambiguity, we drop off the subscript $(i+\frac{1}{2}, j+\frac{1}{2})$ in the notations in the following.

We need to distribute the total residual ${\bf \Phi}$ to the four vertices $\left\{M_{k}\right\}_{k=1,\cdots, 4}$, which is defined in Subsection \ref{SUBSEC:2D_SCALAR} and the corresponding residuals are still denoted by $\left\{{\bf \Phi}^{k}\right\}_{k=1,\cdots, 4}$, where ${\bf \Phi}^{k}\in \mathbb{R}^{m}$. We require that ${\bf \Phi} = \sum^{4}_{k=1}{\bf \Phi}^{k}$ and $|{\bf \Phi}^{k}|/|{\bf \Phi|}, k=1, \cdots, 4$ , which is in a component by component sense and similar to the two-dimensional scalar cases, {\color{orange}should be uniformly bounded for the conservation and the residual property in \cite{Abgrall_JCP2001}. Here we consider a dimension-by-dimension procedure,} coupled with a local characteristic field decomposition. First, we compute an average state $\bar{{\bf u}}$ in $I_{i+\frac{1}{2}, j+\frac{1}{2}}$, using either arithmetic mean or Roe's average \cite{Roe_JCP1981}, and denote $\bar{L}_{x}$ and $\bar{R}_{x}$ as the matrices with left and right eigenvectors $L_{x}$ and $R_{x}$ of ${\bf f}'({\bf u})$ evaluated at the average state $\bar{{\bf u}}$, and $\bar{\lambda}^{k}_{x}$ are the corresponding eigenvalues; {\color{orange}$\bar{L}_{y}$, $\bar{R}_{y}$ and $\bar{\lambda}^{k}_{y}$ are defined similarly but they are associated with $L_{y}$, $R_{y}$ and $\Lambda_{y}$ of ${\bf g}'({\bf u})$, respectively.}

{\color{orange}We now explain how to distribute the total residual within each cell dimension-by-dimension in the upwinding way.}
\begin{description}
\item[Step 1.] {\color{orange}Consider the $y$-direction, and project the residual ${\pmb \Phi}$ to a local characteristic field in the $y$-direction,} we have ${\pmb \Psi}=\bar{L}_{y}{\pmb \Phi}$. {\color{orange}Then the residual ${\pmb \Psi}$ is distributed to the two parts in the $y$-direction,} denoted by ${\pmb \Psi}^{1}$ and ${\pmb \Psi}^{2}$, respectively, and $\pmb \Psi=\pmb \Psi^{1}+\pmb \Psi^{2} $. Residuals ${\pmb \Psi}^{1, 2}$ are defined by
\begin{equation}\label{UW_RD_NWENO2D_SYS_Y}
{\pmb \Psi}^{1}=(\mathrm{I}-\Sigma){\pmb \Psi}, ~~{\pmb \Psi}^{2}=\Sigma{\pmb \Psi},
\end{equation}
where $\mathrm{I}$ is the identity matrix, $\Sigma$ is a diagonal matrix with the $k$-th diagonal component given by
\begin{equation*}
\Sigma_{kk}=\left\{\begin{array}{ll}
1 & \mbox{if}~\bar{\lambda}^{k}_{y} \geq \delta, \\
0 & \mbox{if}~\bar{\lambda}^{k}_{y} \leq -\delta, \\
r(\bar{\lambda}^{k}_{y}, \delta) & \mbox{otherwise}.
\end{array}\right.
\end{equation*}
The function $r(\cdot, \cdot)$ is given in \eqref{ROE_FUN}, and $\delta$ is also chosen accordingly in the problem. And then project residuals ${\pmb \Psi}^{1,2}$ back to the physical space, we obtain residuals $\hat{{\pmb \Phi}}^{1, 2}$, namely
\begin{equation}\label{UW_HWENO2D_SYS_PB_Y}
\hat{{\pmb \Phi}}^{1}=\bar{R}_{y}{\pmb \Psi}^{1}, ~~\hat{{\pmb \Phi}}^{2}=\bar{R}_{y}{\pmb \Psi}^{2}.
\end{equation}

\item[Step 2.] Consider the $x$-direction, and we would distribute the two parts $\hat{{\pmb \Phi}}^{1, 2}$ in the $x$-direction. First we need to project residuals $\hat{{\pmb \Phi}}^{1, 2}$ to the characteristic fields in the $x$-direction, namely
\begin{equation*}
\Pi^{1}=\bar{L}_{x}\hat{{\pmb \Phi}}^{1}, ~~\Pi^{2}=\bar{L}_{x}\hat{{\pmb \Phi}}^{2}.
\end{equation*}
Then distribute residuals $\Pi^{1,2}$ in the $x$-characteristic fields. {\color{orange}According to the  upwinding principle and the residual property,} we have
\begin{equation}\label{UW_RD_HWENO2D_SYS_X}
\bar{{\pmb \Psi}}^{1} =(\mathrm{I}-\Gamma)\Pi^{1},~\bar{{\pmb \Psi}}^{2} =\Gamma\Pi^{1},~\bar{{\pmb \Psi}}^{3} =(\mathrm{I}-\Gamma)\Pi^{2},~\bar{{\pmb \Psi}}^{4} =\Gamma\Pi^{2},
\end{equation}
where $\mathrm{I}$ is the identity matrix, and $\Gamma$ is a diagonal matrix with the $k$-th diagonal component given by
\begin{equation*}
\Gamma_{kk}=\left\{\begin{array}{ll}
1 & \mbox{if}~\bar{\lambda}^{k}_{x} \geq \delta, \\
0 & \mbox{if}~\bar{\lambda}^{k}_{x} \leq -\delta, \\
r(\bar{\lambda}^{k}_{x}, \delta) & \mbox{otherwise}.
\end{array}\right.
\end{equation*}
The function $r(\cdot, \cdot)$ is given as in \eqref{ROE_FUN}, and $\delta$ is chosen accordingly in the problem.

\item[Step 3.] {\color{orange}Project distributed residuals $\left\{\bar{{\pmb \Psi}}^{k}\right\}_{k=1, \cdots, 4}$ back to the physical space,} namely
\begin{equation}\label{UW_HWENO2D_SYS_PB_X}
\tilde{{\pmb \Phi}}^{k}=\bar{R}_{x}\bar{{\pmb \Psi}}^{k}, ~~k=1, \cdots, 4.
\end{equation}
{\color{red}As in two-dimensional scalar cases, we need to add an additional dissipation residual ${\pmb \Phi}^{k}_{\mbox{diss}}, k=1, \cdots, 4$ to each of $\pmb \Phi^{k}, k=1, \cdots, 4$}, and its definition is defined as in \eqref{DEF:UW_HWENO2D_RESIDUAL_DISS}. Thus, we get the way to distribute the total residuals
\begin{equation}\label{UW_HWENO2D_RES_FINAL_SYS}
{\pmb \Phi}^{k}=\tilde{{\pmb \Phi}}^{k}+{\pmb \Phi}^{k}_{\mbox{diss}}, ~~k=1, \cdots, 4.
\end{equation}
\end{description}
{\color{orange}The point value ${\bf u}^{\text{new}}_{i, j}$ is updated through sending the distributed residuals to the point $(x_{i}, y_{j})$,} as in the pseudo time-marching scheme, which can be written as the semi-discrete systems \eqref{SCHEMA:UW_HWENO2D_SCALAR}. And we again use a forward Euler scheme for the pseudo time discretization in our numerical experiments until the steady state is reached. {\color{orange}And for ${\bf v}^{\text{new}}_{i, j}$ and ${\bf w}^{\text{new}}_{i, j}$, they are again updated by the fourth order HWENO reconstruction component by component in the $x$-direction and in the $y$-direction, respectively. ${\bf z}^{\text{new}}_{i, j}$ is updated by the fourth order HWENO reconstructions in a dimension-by-dimension way.}

{\color{orange}
\underline{\bf Procedure \rom{4}. RD finite difference HWENO method for steady state problem in 2D}
\begin{enumerate}
\item Compute the total residuals within the cell component by component as defined in \eqref{DEF:RES2D}, in which the integrals of fluxes are approximated by the sixth order HWENO integration as in the one-dimensional case, the integral of the source term is approximated in a dimension-by-dimension via the sixth order HWENO integration.

\item Project the total residual to a local characteristic field in the $y$-direction, and then distribute the residual into two parts in the $y$-direction in an upwinding way, as in \eqref{UW_RD_NWENO2D_SYS_Y}, and then project the distributed residuals back to the physical space, as in \eqref{UW_HWENO2D_SYS_PB_Y}.

\item Project the two parts of the residuals to a local characteristic field in the $x$-direction, and distribute them to the four vertices of the cell, as in \eqref{UW_RD_HWENO2D_SYS_X}, according to the upwinding principle. And then project the distributed residuals back to the physical space, as in \eqref{UW_HWENO2D_SYS_PB_X}.

\item Revise the distributed residuals of the four vertices by adding an additional dissipation residuals, as shown in \eqref{UW_HWENO2D_RES_FINAL_SYS}.

\item Update the point value ${\bf u}^{\text{new}}_{i, j}$ through sending the residuals in the physical space and forward in pseudo time \eqref{SCHEMA:UW_HWENO2D_SCALAR} by a forward Euler time discretization until the steady state is reached.

\item Update the spatial derivatives of the point value ${\bf v}^{\text{new}}_{i,j}$  and ${\bf w}^{\text{new}}_{i, j}$ by the fourth order HWENO reconstruction in the $x$-direction and in the $y$-direction, respectively, as in Subsection \ref{SUBSEC:1D_SCALAR}. The second cross derivative ${\bf z}^{\text{new}}_{i, j}$ is approximated by the fourth order HWENO reconstruction in a dimension-by-dimension way.
\end{enumerate}
}

\section{Numerical results}
In this section, {\color{orange}we present the numerical results of the proposed residual distribution conservative finite difference HWENO method for steady state conservation laws with source terms in scalar and system problems in one and two dimensions. Pseudo time discretization towards steady state is by the forward Euler method in  all numerical simulations.} CFL number is taken to be 0.6 in one-dimensional cases and 0.2 in two-dimensional cases. In Subsection \ref{SUBSEC:1D_SCAL}, for one-dimensional problems, {\color{blue}the parameter $\delta$ in \eqref{ROE_FUN} for the Roe's entropy correction is taken as $10^{-15}$. In Subsection \ref{SUBSEC:2D_SCAL}, for two-dimensional problems, $\delta$ is taken as $10^{-15}$ for scalar cases and $0.1$ for system cases, unless otherwise stated.}

All the spatial discretizations in our numerical results are uniform and all numerical steady state is obtained with $L^{1}$ residue reduced to the round-off level.

\subsection{The one-dimensional problems}\label{SUBSEC:1D_SCAL}
\begin{exa}\label{EG:UW_HWENO_BUR1D_SM}
{\color{red}We solve the steady state solution of the one-dimensional  Burgers' equation with a source term:}
\begin{equation}\label{EQ:UW_HWENO_BUR1D_SM}
\left(\frac{u^{2}}{2}\right)_{x}=\sin x \cos x
\end{equation}
with the initial condition
\begin{equation}\label{EQ:UW_HWENO_BUR1D_SM_INIT}
u_{0}(x) = \beta \sin x,
\end{equation}
and the boundary condition $u(0) = u(\pi) = 0$. This problem was studied in \cite{S.A.Gottlieb_ANM1986} as an example of multiple steady state solutions for characteristic initial value problems. The steady state solution to this problem depends on the value of $\beta$: if $-1<\beta<1$, a shock will form within the domain $\left[0,\pi\right]$; otherwise, the solution will be smooth at first, followed by a shock forming at the boundary $x = \pi ~(\beta \geq 1)$ or $x = 0~(\beta \leq -1)$, and later converge to a smooth steady state $u(x, \infty) = \sin x~(\beta \geq 1)$ or $u(x, \infty) = -\sin x~(\beta \leq -1)$, respectively. In order to test the order of accuracy, we take $\beta = 2$ to have a smooth stationary solution. {\color{orange}From Table \ref{TAB:UW_HWENO_BUR1D_SM},} we can clearly see that the sixth order accuracy is reached. {\color{orange}In Figure \ref{PIC:UW_HWENO-WENO_ERROR_BUR1D},} {\color{blue}we can observe that the $L^{1}$ error of the proposed scheme combined with  the novel HWENO scheme  is smaller than that of the RD scheme combined with the WENO scheme and very close to that of the RD scheme in the traditional HWENO framework at the same grid.}
\begin{table}[!htbp]
\centering
\caption{{\color{orange}Errors and numerical orders of accuracy for the sixth order RD finite difference HWENO scheme in Example \ref{EG:UW_HWENO_BUR1D_SM}.}}
\vspace{5mm}
\begin{tabular}{lllll}
   \hline
    $N$ & $L^{1}$ \mbox{error} & \mbox{Order} & $L^{\infty}$ \mbox{error} & \mbox{Order} \\ \hline
    20 &     4.69E-07 &          &     8.30E-07 &          \\
    40 &     4.31E-09 &     6.76 &     8.03E-09 &     6.69 \\
    80 &     3.96E-11 &     6.77 &     7.55E-11 &     6.73 \\
    160 &     3.78E-13 &     6.71 &     7.07E-13 &     6.74 \\
    320 &     3.93E-15 &     6.59 &     6.87E-15 &     6.69 \\
    640 &     4.52E-17 &     6.44 &     7.42E-17 &     6.53 \\
  \hline
\end{tabular}
\label{TAB:UW_HWENO_BUR1D_SM}
\end{table}

\begin{figure}[!htbp]
	\centering
	\includegraphics[width=0.45\textwidth]{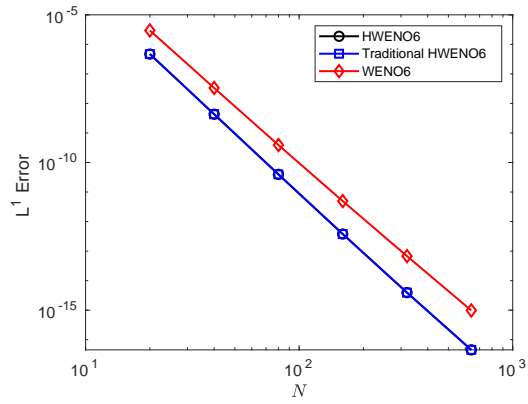}	
	\caption{{\color{blue}$L^{1}$ Error for HWENO6, traditional HWENO6 and  WENO6 schemes in Example \ref{EG:UW_HWENO_BUR1D_SM}.}}
	\label{PIC:UW_HWENO-WENO_ERROR_BUR1D}
\end{figure}
\end{exa}

\begin{exa}\label{EG:UW_HWENO_BUR1D_NSM}
{\color{orange}We consider the same problem as Example \ref{EG:UW_HWENO_BUR1D_SM},} but here take $\beta = 0.5$ in the initial condition \eqref{EQ:UW_HWENO_BUR1D_SM_INIT}. As mentioned in the previous example, when $-1 < \beta < 1$, a shock will form within the domain, which separates two branches ($\sin x$ and $-\sin x$) of the steady state. The location of the shock is determined by the parameter $\beta$ through conservation of mass ($\int^{\pi}_{0} u \,dx = 2\beta$), and can be derived to be $\pi - \arcsin \sqrt{1-\beta^{2}}$. For the case $\beta = 0.5$, the shock location is approximately 2.0944. {\color{blue}The numerical solution on the uniform meshes is shown in Figure \ref{PIC:UW_HWENO_BUR1D_NSM} (left). We can see that the numerical shock is at the correct location and is resolved well. From Figure \ref{PIC:UW_HWENO_BUR1D_NSM} (right), CPU time of the novel HWENO scheme is very close to that of the traditional HWENO scheme, as the grid increases, but the WENO scheme takes less time to reach steady state.}	
%
\begin{figure}[!htbp]
  \centering%
    \includegraphics[width=0.45\textwidth]{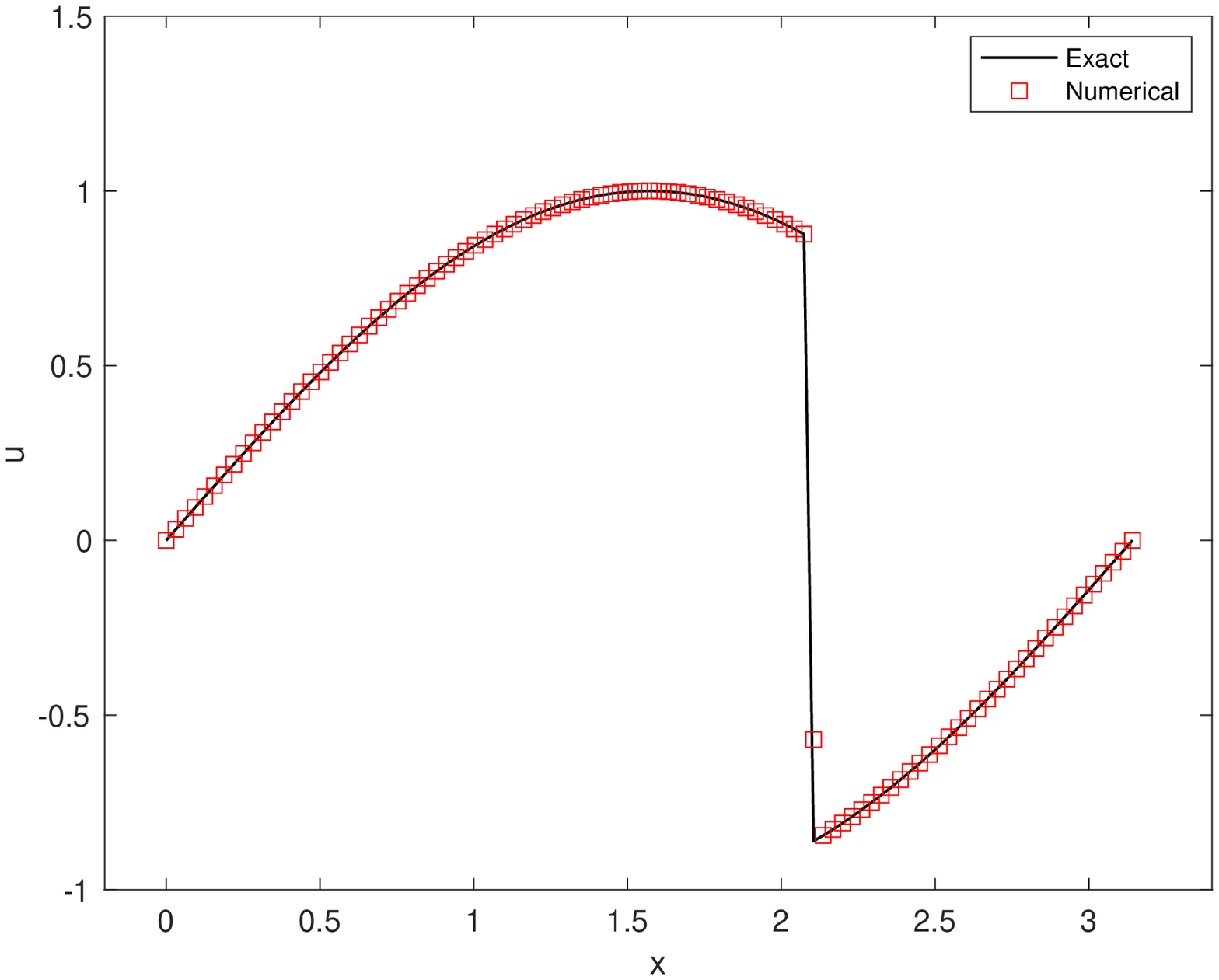}
	\includegraphics[width=0.45\textwidth]{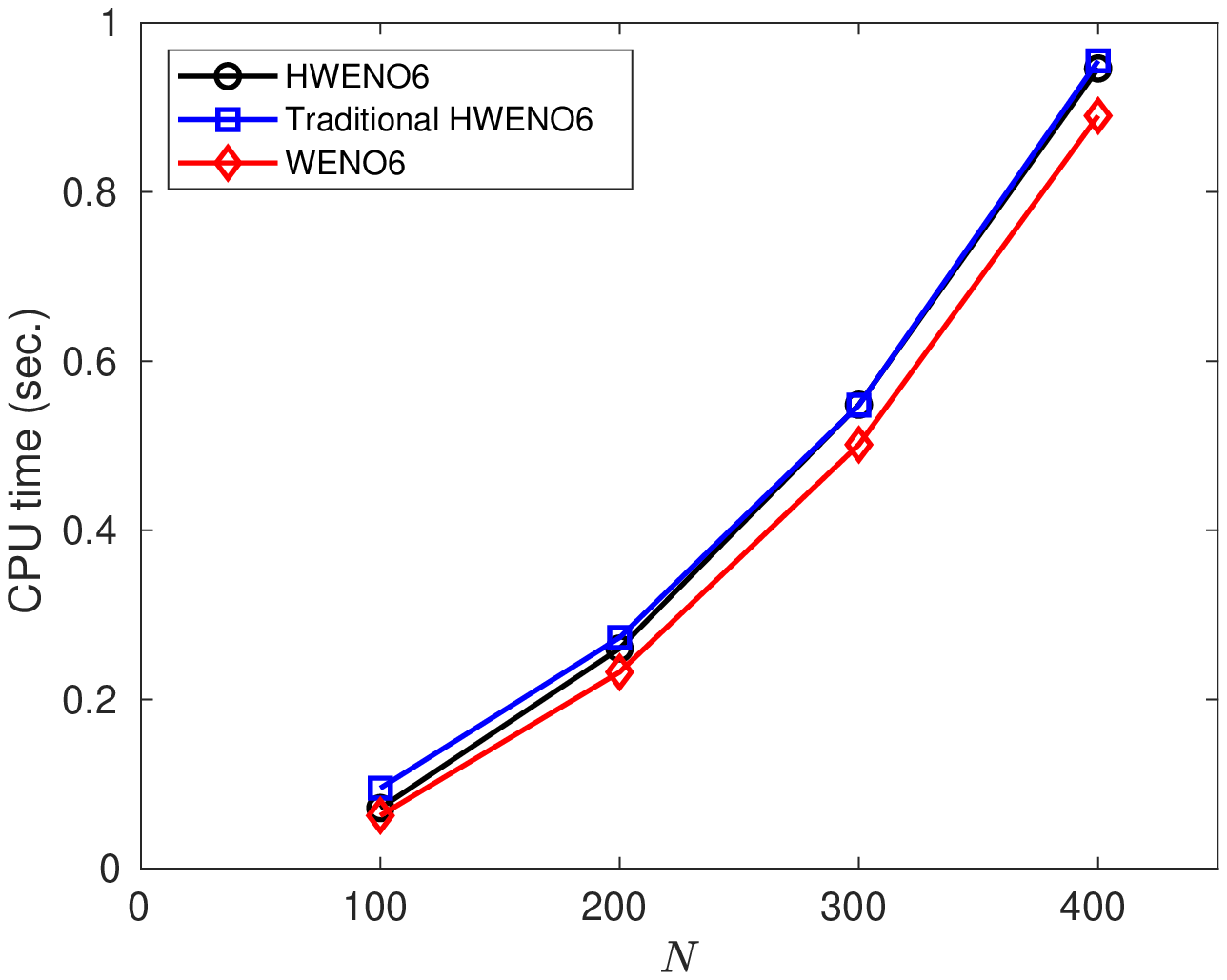}
    \caption{{\color{blue}Example \ref{EG:UW_HWENO_BUR1D_NSM}. Left: the numerical solution vs. the exact solution with 100 cells; right: CPU time for HWENO6, traditional HWENO6 and WENO6 schemes.}}
  \label{PIC:UW_HWENO_BUR1D_NSM}
\end{figure}
\end{exa}

\begin{exa}\label{EG:UW_HWENO_BUR1D_NSM2}
{\color{red} We consider the steady state solutions of the  Burgers' equation with a different source term}, which depends on the solution itself:
\begin{equation}
\left( \frac{u^{2}}{2}\right)_{x} = -\pi \cos(\pi x)u, ~~ x\in \left[0, 1\right]
\end{equation}
equipped with the boundary conditions $u(0) = 1$ and $u(1) = -0.1$. This problem has two steady state solutions with shocks
\[
u(x) = \begin{cases}
u^{+} = 1 - \sin(\pi x) & \text{if $0 \leq x < x_{s}$}, \\
u^{-} = -0.1 - \sin(\pi x) & \text{if $x_{s} \leq x \leq 1$},
\end{cases}
\]
where $x_{s} = 0.1486$ or $x_{s} = 0.8514$. Both solutions satisfy the Rankine-Hugoniot jump condition and the entropy conditions, but only the one with the shock at 0.1486 is stable for a small perturbation. This problem was studied in \cite{Embid.Goodman.Majda_SIAM_SSC1984} as an example of multiple steady states for one-dimensional transonic flows. This case is tested to demonstrate that starting with a reasonable perturbation of the stable steady state, the numerical solution converges to the stable one.

The initial condition is given by
\[
u_{0}(x) = \begin{cases}
\hfil 1 & \text{if $0 \leq x < 0.5$},\\
\hfil -0.1  & \text{if $0.5 \leq x \leq 1$},
\end{cases}
\]
where the initial jump is located in the middle of the position of the shocks in the two admissible steady state solution. {\color{orange}The numerical result and the exact solution are displayed in Figure \ref{PIC:UW_HWENO_BUR1D_NSM2}.} We can see the shock location and good resolution of the shock. {\color{orange}And we also obtain the sixth order accuracy of our scheme in a smooth region $[0.5, 1]$ from Table \ref{TAB:UW_HWENO_BUR1D_NSM2}.} Note that in this case the source term is dependent to the numerical solution.

\begin{figure}[!htbp]
	\centering%
	\includegraphics[width=0.45\textwidth]{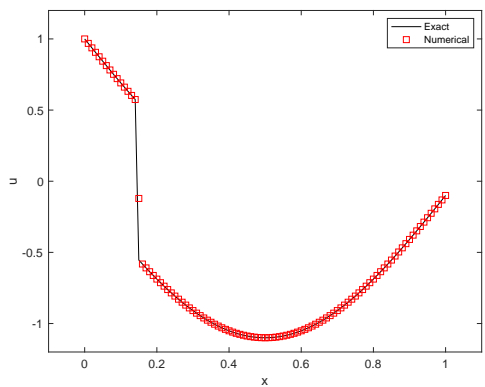}
	
	\caption{{\color{orange}Example \ref{EG:UW_HWENO_BUR1D_NSM2}. The numerical solution vs. the exact solution with 100 cells.}}
	\label{PIC:UW_HWENO_BUR1D_NSM2}
\end{figure}
%

\begin{table}[!htbp]
	\centering
	\caption{{\color{orange}Errors and numerical orders of accuracy for the sixth order RD finite difference HWENO scheme in Example \ref{EG:UW_HWENO_BUR1D_NSM2} at $[0.5, 1]$.}}
	\vspace{5mm}
	\begin{tabular}{lllll}
		\hline
		$N$ & $L^{1}$ \mbox{error} & \mbox{Order} & $L^{\infty}$ \mbox{error} & \mbox{Order} \\ \hline
		20 &     4.56E-07 &          &     7.30E-07 &          \\
		40 &     6.32E-09 &     6.17 &     1.50E-08 &     5.61 \\
		80 &     1.63E-10 &     5.28 &     2.78E-10 &     5.76 \\
		160 &     3.08E-12 &     5.72 &     4.78E-12 &     5.86 \\
		320 &     5.25E-14 &     5.88 &     7.85E-14 &     5.93 \\
		640 &     8.54E-16 &     5.94 &     1.26E-15 &     5.96 \\
		\hline
	\end{tabular}
	\label{TAB:UW_HWENO_BUR1D_NSM2}
\end{table}
\end{exa}

{\color{blue}
	\begin{rem}
		In \cite{renhigh}, the author has shown that the novel HWENO framework has the advantages of less storage and lower cost than the traditional HWENO for solving static Hamilton-Jacobi equations, but this advantage is not seen here in our simulations, see Example \ref{EG:UW_HWENO_BUR1D_NSM}. The reason is that the traditional HWENO framework does not need to carry out HWENO reconstruction when updating $v^{\text{new}}$ in one-dimensional case, only the residual distribution of the auxiliary equation is required. However, $v^{\text{new}}$ is updated by HWENO reconstruction in this paper, so that it takes longer time, even though the method does not introduce the auxiliary equation. 
	\end{rem}
}

\begin{exa}\label{EG:UW_HWENO_SWATER1D}
We solve the steady state solutions of the one-dimensional shallow water equation

\begin{equation}
\left(\begin{array}{c}
hu \\
hu^{2}+\frac{1}{2}gh^{2}

\end{array}\right)_{x}
=\left(\begin{array}{c}
0 \\
-ghb_{x}
\end{array}\right),
\end{equation}
where $h$ denotes the water height, $u$ is the velocity of the fluid, $b(x)$ represents the bottom topography and $g$ is the gravitational constant.

Starting from a stationary initial condition, which itself is a steady state solution, we can check the order of accuracy. The smooth bottom topography is given by
\begin{equation*}
b(x) = 5\exp^{-\frac{2}{5}(x - 5)^{2}}, ~~x\in [0, 10].
\end{equation*}
The initial condition is the stationary solution
{\color{orange}
\begin{equation*}
h + b = 10, ~~hu = 0,
\end{equation*}}
and the exact steady state solution is imposed as the boundary condition.

We test our scheme on uniform meshes. {\color{orange}The numerical accuracies are shown in Table \ref{TAB:UW_HWENO_SWATER1D}. We can clearly see the orders of accuracies and errors for the water height $h$.} {\color{blue} Figure \ref{PIC:UW_HWENO-WENO_ERROR_	SWATER} shows that the comparisons of $L^{1}$ error among the novel HWENO, traditional HWENO  and  WENO schemes. We can clearly observe that the $L^{1}$ error of the proposed scheme combined with  the novel HWENO scheme  is smaller than that of the RD scheme combined with the WENO scheme and very close to that of the RD scheme in the traditional HWENO framework at the same grid.}
\begin{table}[!htbp]
\centering
\caption{{\color{orange}Errors and numerical orders of accuracy for the water height $h$ of the sixth order RD finite difference HWENO scheme in Example \ref{EG:UW_HWENO_SWATER1D}.}}
\vspace{5mm}
\begin{tabular}{lllll}
  \hline
    $N$ & $L^{1}$ \mbox{error} & \mbox{Order} & $L^{\infty}$ \mbox{error} & \mbox{Order} \\ \hline
     20 &     7.22E-04 &          &     2.78E-03 &          \\
     40 &     9.07E-06 &     6.31 &     3.63E-05 &     6.26 \\
     80 &     1.03E-07 &     6.47 &     4.38E-07 &     6.37 \\
    160 &     1.07E-09 &     6.58 &     4.92E-09 &     6.48 \\
    320 &     1.12E-11 &     6.58 &     5.58E-11 &     6.46 \\
    640 &     1.26E-13 &     6.48 &     6.80E-13 &     6.36 \\
  \hline
\end{tabular}
\label{TAB:UW_HWENO_SWATER1D}
\end{table}

\begin{figure}[!htbp]
	\centering
	\includegraphics[width=0.45\textwidth]{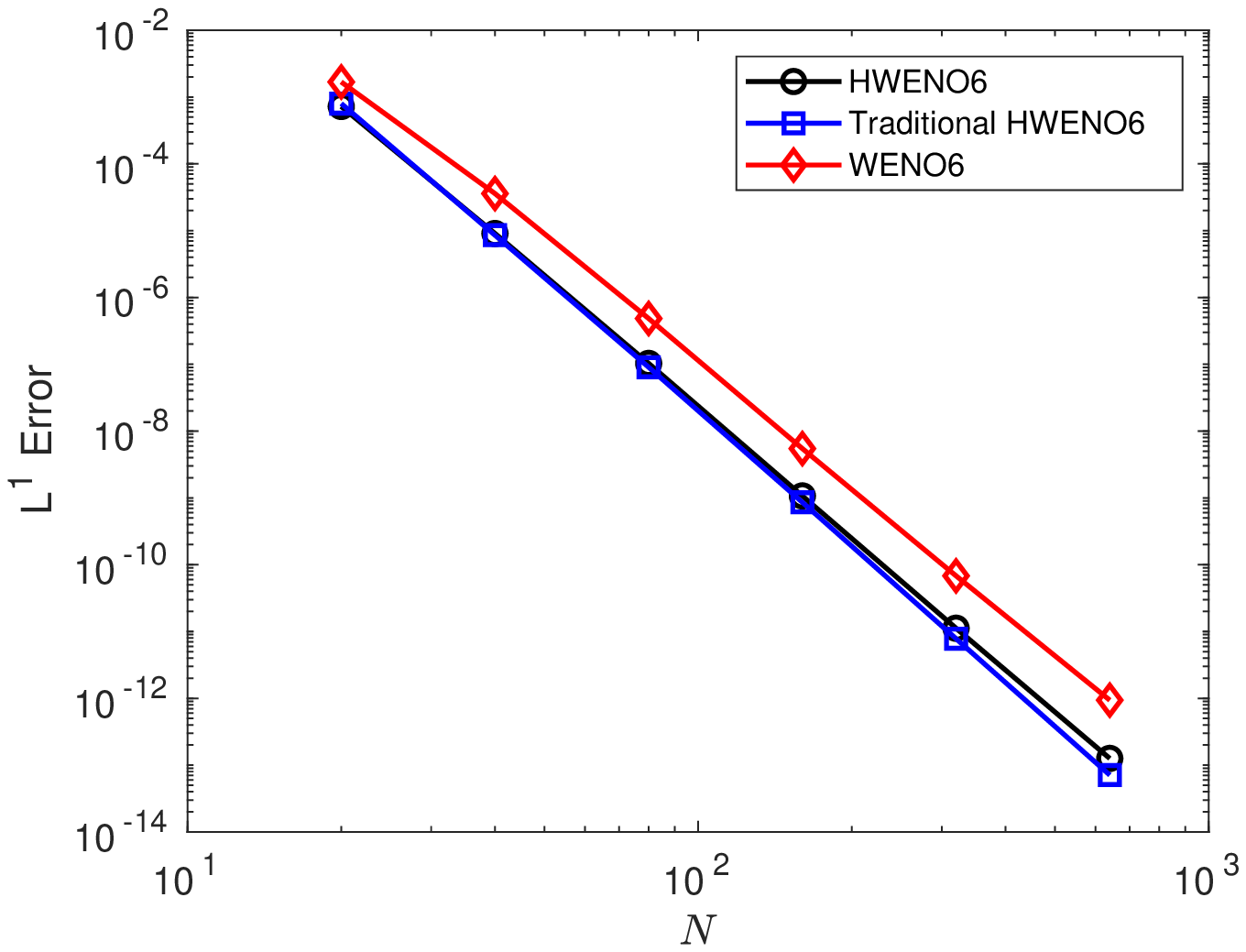}	
	\caption{{\color{blue}$L^{1}$ Error for HWENO6, traditional HWENO6  and WENO6 schemes in Example \ref{EG:UW_HWENO_SWATER1D}.}}
	\label{PIC:UW_HWENO-WENO_ERROR_	SWATER}
\end{figure}
\end{exa}

\begin{exa}\label{EG:UW_HWENO_NOZZLE1D}
We test our scheme on the steady state solution of the one-dimensional nozzle flow problem
\begin{equation}
\left(\begin{array}{c}
\rho u \\
\rho u^{2} + p \\
u(E + p)
\end{array}\right)_{x}
=-\frac{A'(x)}{A(x)}\left(\begin{array}{c}
\rho u \\
\rho^{2}u^{2}/\rho \\
u(E + p)
\end{array}\right), ~~ x\in \left[0, 1\right],
\end{equation}
where $\rho$ is the density, $u$ is the velocity of the fluid, $E$ is the total energy, $\gamma=1.4$ is the gas constant, $p = (\gamma -1)(E - \frac{1}{2}\rho u^{2})$ is the pressure and $A(x)$ represents the area of the cross-section of the nozzle.

We start with an isentropic initial condition, with a shock at $x = 0.5$. The density $\rho$ and pressure $p$ at $-\infty$ are 1, and the inlet Mach number at $x = 0$ is 0.8. The outlet Mach number at $x = 1$ is 1.8, with linear Mach number distribution before and after the shock. The area of the cross-section $A(x)$ is then determined by the relation
\begin{equation*}
A(x)f(\text{Mach number at x}) = \text{constant}, ~~ \forall x\in [0, 1],
\end{equation*}
where
\begin{equation*}
f(w) = \frac{w}{(1+\delta_{0} w^{2})^{p_{0}}},~~\delta_{0}=\frac{1}{2}(\gamma -1),~~p_{0} = \frac{1}{2}\cdot\frac{\gamma + 1}{\gamma - 1}.
\end{equation*}

{\color{orange}From Figure \ref{PIC:UW_HWENO_NOZZLE1D}, we can clearly see that the shock is resolved well.}
\begin{figure}[!htbp]
  \centering
  \mbox{
    \includegraphics[width=0.45\textwidth]{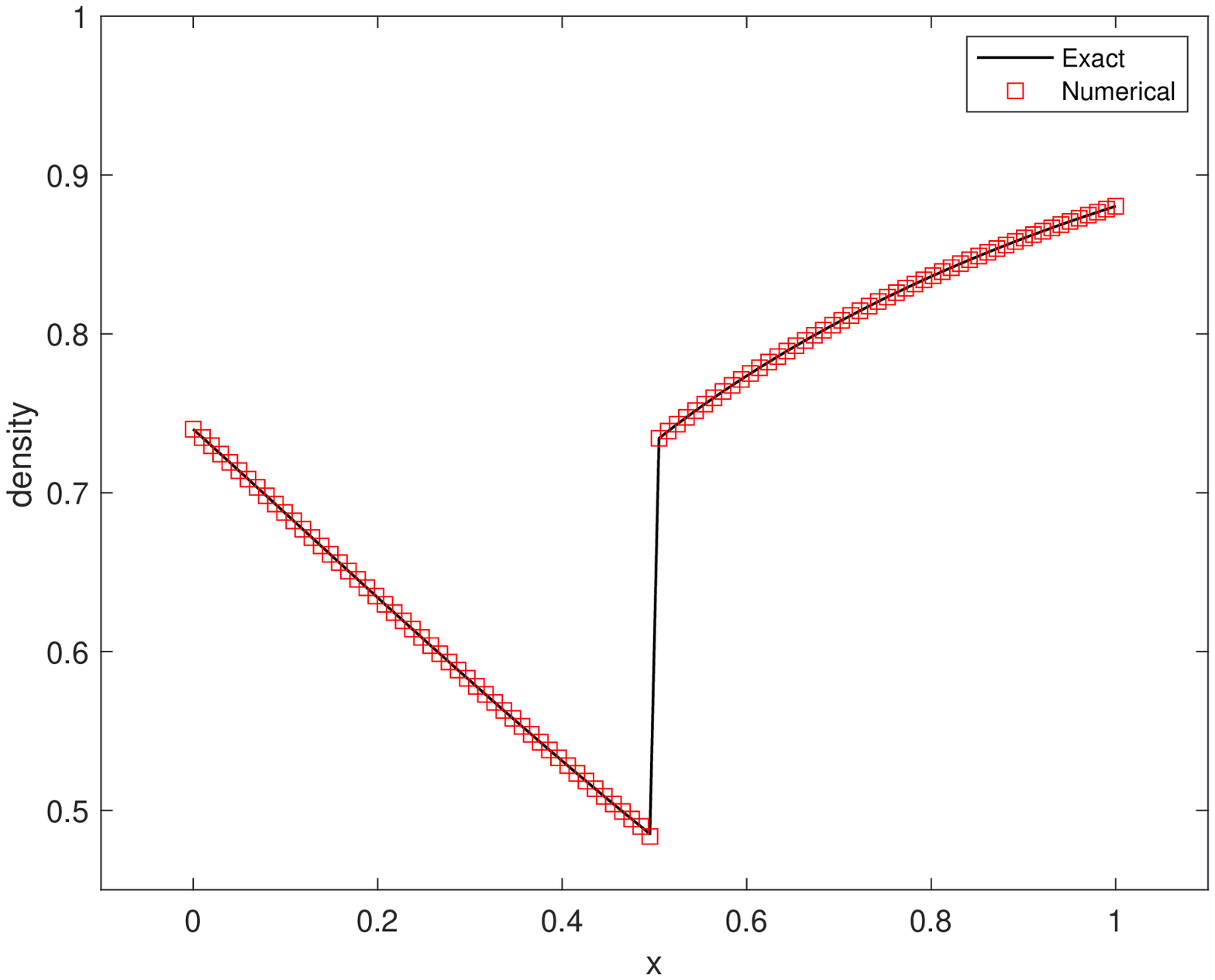}
  	\quad
    \includegraphics[width=0.45\textwidth]{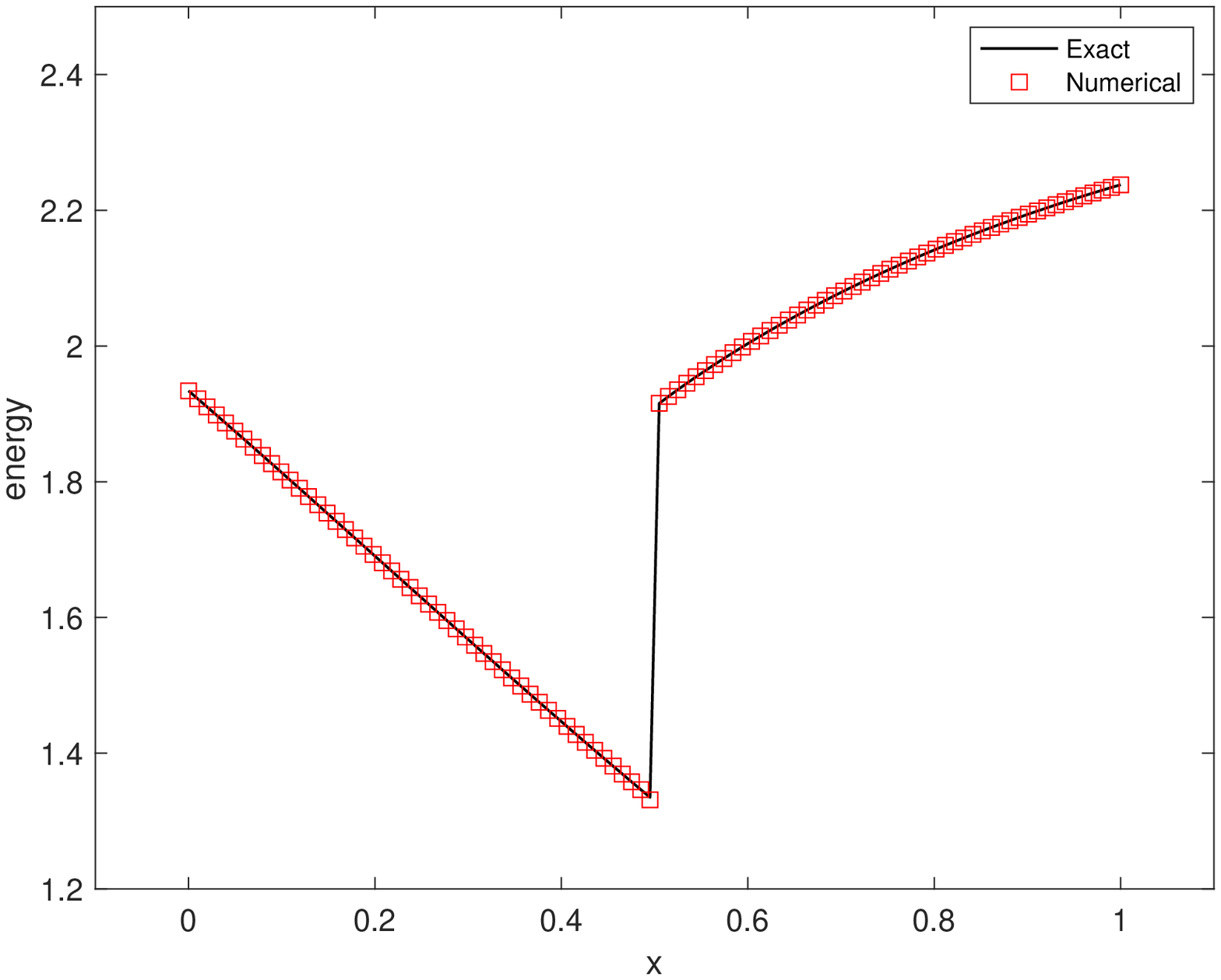}}
  \mbox{
    \includegraphics[width=0.45\textwidth]{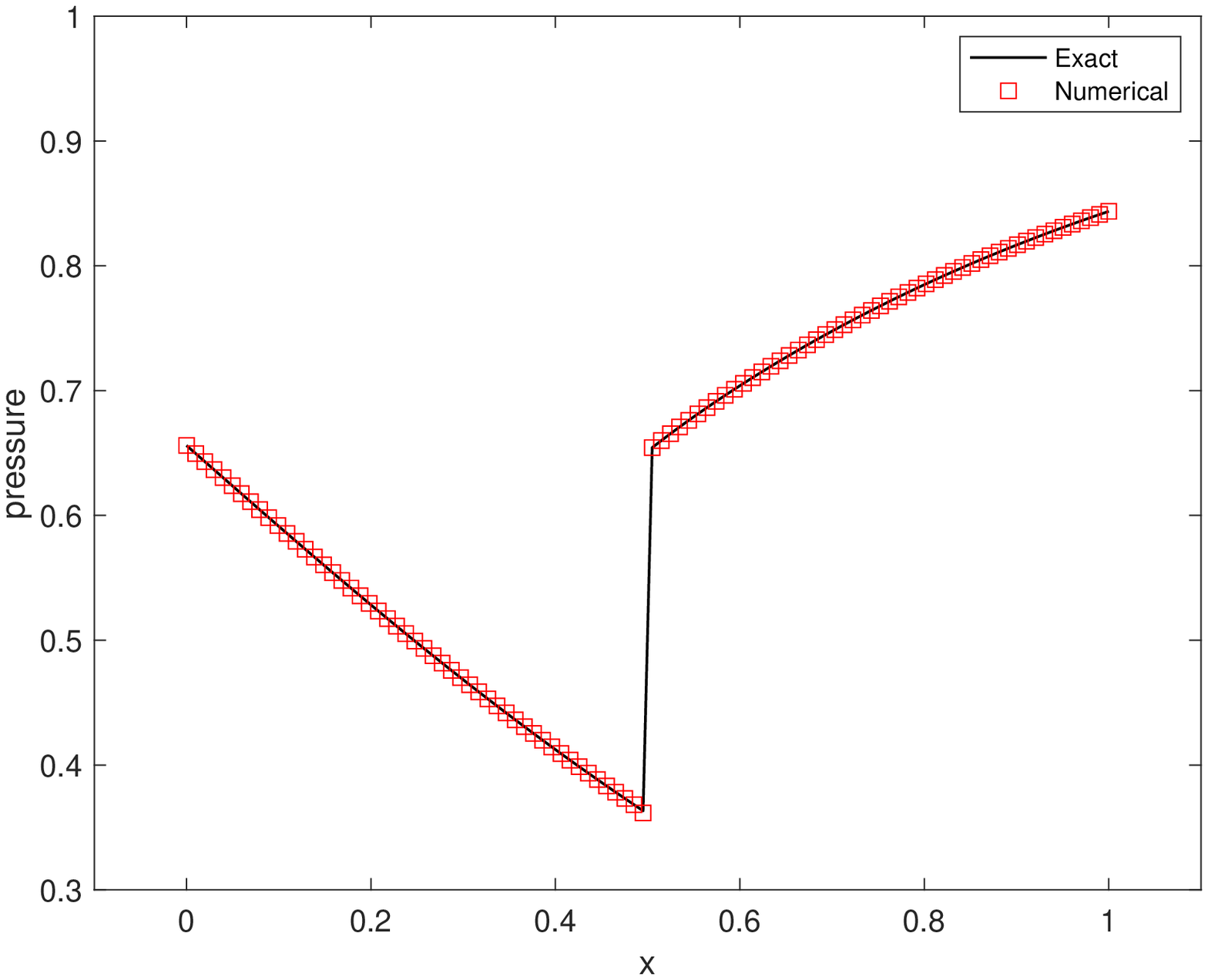}
  	\quad
    \includegraphics[width=0.45\textwidth]{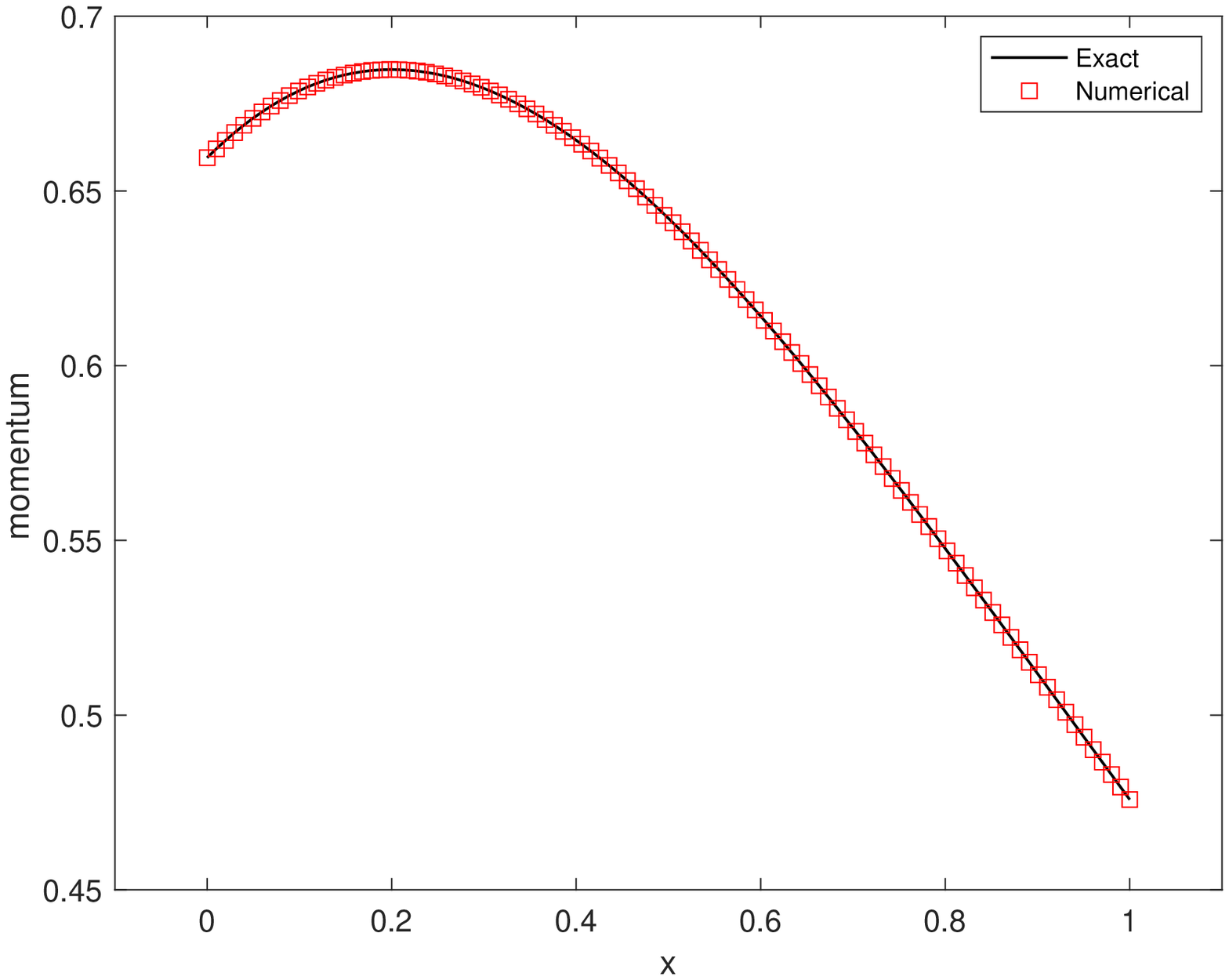}}

  \caption{{\color{orange}Nozzle flow problem with 101 cells in Example \ref{EG:UW_HWENO_NOZZLE1D}. Top left: density; top right: total energy; bottom left: pressure; bottom right: momentum.}}
  \label{PIC:UW_HWENO_NOZZLE1D}
\end{figure}

\end{exa}

\subsection{The two-dimensional problems}\label{SUBSEC:2D_SCAL}

\begin{exa}\label{EG:UW_HWENO_BUR2D_SM}
{\color{red} We solve the steady state problem of two-dimensional Burgers' equation with a source term}
\begin{equation}
\left(\frac{1}{\sqrt{2}}\frac{u^{2}}{2}\right)_{x} + \left(\frac{1}{\sqrt{2}}\frac{u^{2}}{2}\right)_{y} = \sin\left(\frac{x+y}{\sqrt{2}}\right) \cos\left(\frac{x+y}{\sqrt{2}}\right),
\end{equation}
where $(x, y) \in \left[0, \frac{\pi}{\sqrt{2}}\right]\times \left[0, \frac{\pi}{\sqrt{2}}\right]$ with the initial condition given by
\begin{equation}
u_{0}(x, y) = \beta \sin\left(\frac{x+y}{\sqrt{2}}\right).
\end{equation}
{\color{orange}This is the one-dimensional problem studied in Example \ref{EG:UW_HWENO_BUR1D_SM} along the northeast-southwest diagonal.} Since our grids are not aligned with the diagonal, this is a truly two-dimensional test case. Here we take the boundary conditions to be the exact solution of the steady state problem.

For this example, we take $\beta = 1.2$, which gives a smooth steady state solution $u(x, y) = \sin\left(\frac{x+y}{\sqrt{2}}\right)$. The errors and numerical orders are shown in the Table \ref{TAB:UW_HWENO_BUR2D_SM}. {\color{red}It can be clearly seen that the sixth order accuracy is achieved in the $L^{1}$ error case principally.} {\color{blue}But, when the grid increases, we can also achieve the sixth order accuracy in the $L^{\infty}$ sense.}
\begin{table}[!htbp]
\centering
\caption{{\color{blue}Errors and numerical orders of accuracy for the sixth order RD finite difference HWENO scheme in Example \ref{EG:UW_HWENO_BUR2D_SM} with $N \times N $ cells.}}
\vspace{5mm}
\begin{tabular}{llllllll}
  \hline
    $N\times N$ & $L^{1}$ \mbox{error} & \mbox{Order} & $L^{2}$ \mbox{error} & \mbox{Order} & $L^{\infty}$ \mbox{error} & \mbox{Order} \\ \hline
$    10\times    10$ &     8.37E-07 &          &     9.82E-07 &          &     2.28E-06 &          \\
$    20\times    20$ &     8.76E-09 &     6.58 &     1.10E-08 &     6.48 &     5.31E-08 &     5.42 \\
$    40\times    40$ &     9.72E-11 &     6.50 &     1.40E-10 &     6.30 &     1.52E-09 &     5.13 \\
$    80\times    80$ &     1.13E-12 &     6.42 &     2.00E-12 &     6.12 &     4.63E-11 &     5.04 \\
$   160\times   160$ &     1.33E-14 &     6.41 &     3.03E-14 &     6.05 &     1.37E-12 &     5.08 \\
$   320\times   320$ &     1.73E-16 &     6.26 &     3.53E-16 &     6.42 &     1.98E-14 &     6.11 \\
  \hline
\end{tabular}
\label{TAB:UW_HWENO_BUR2D_SM}
\end{table}
\end{exa}

{\color{blue}
	\begin{rem}	
		We want to emphasize that here we cannot compare the numerical results of the two HWENO frameworks in two-dimensional problems in our paper. In \cite{mythesis}, the author only developed one-dimensional RD scheme in the traditional HWENO framework, because it is still difficult in residual distribution for the auxiliary equations and will be explored in future. In principle, for the two-dimensional problem, three auxiliary equations need to be introduced under the traditional HWENO framework, including four integral terms. Therefore, at least four times sixth order HWENO integration procedures are required to solve the auxiliary equations to obtain $v$. However, the method in this paper only requires three times fourth order HWENO reconstructions.  Therefore, there is less storage and costs under the novel HWENO framework.
	\end{rem}
}

\begin{exa}\label{EG:UW_HWENO_BUR2D_NSM}
We consider the steady state solution of the following problem:
\begin{equation}
\left(\frac{1}{\sqrt{2}}\frac{u^{2}}{2}\right)_{x} + \left(\frac{1}{\sqrt{2}}\frac{u^{2}}{2}\right)_{y} = -\pi \cos(\pi \frac{x+y}{\sqrt{2}})u,
\end{equation}
where $(x, y)\in \left[0, \frac{1}{\sqrt{2}}\right]\times
\left[0, \frac{1}{\sqrt{2}}\right]$. {\color{orange}This is the one-dimensional problem in  Example \ref{EG:UW_HWENO_BUR1D_NSM} along the northeast-southwest diagonal line.} Inflow boundary conditions are given by the exact solution of the steady state problem. Again, since our grids are not aligned with the diagonal line, this is a truly two-dimensional test case. As before, this problem has two steady state solutions with shocks
\begin{equation*}
u(x , y) =\begin{cases}
1- \sin\left(\pi \frac{x+y}{\sqrt{2}}\right)  & \text{if $0 \leq \frac{x+y}{\sqrt{2}} < x_{s}$}, \\
-0.1 - \sin\left(\pi \frac{x+y}{\sqrt{2}}\right) & \text{if $x_{s} \leq \frac{x+y}{\sqrt{2}} \leq 1$},
\end{cases}
\end{equation*}
where $x_s = 0.1486$ or $x_{s} = 0.8514$. Both solutions satisfy the Rankine-Hugoniot jump condition and the entropy conditions, but only the one with the shock at $\frac{x+y}{\sqrt{2}}= 0.1486$ is stable for a small perturbation.

The initial condition is given by
\[
u_{0}(x, y) = \begin{cases}
\hfil 1 & \text{if $0\leq \frac{x+y}{\sqrt{2}} < 0.5$}, \\
\hfil -0.1 & \text{if $0.5 \leq \frac{x+y}{\sqrt{2}} \leq 1$},
\end{cases}
\]
where the initial jump is located in the middle of the positions of the shocks in the two admissible steady state solutions. From Figure \ref{PIC:UW_HWENO_BUR2D_NSM}, we can see the correct shock location and a good resolution of the solution. {\color{orange}The coefficient $\sigma$ for the dissipation \eqref{DEF:UW_HWENO2D_RESIDUAL_DISS} is taken as $4$.}
\begin{figure}[!htbp]
  \centering%
  \mbox{
    \includegraphics[width=0.45\textwidth]{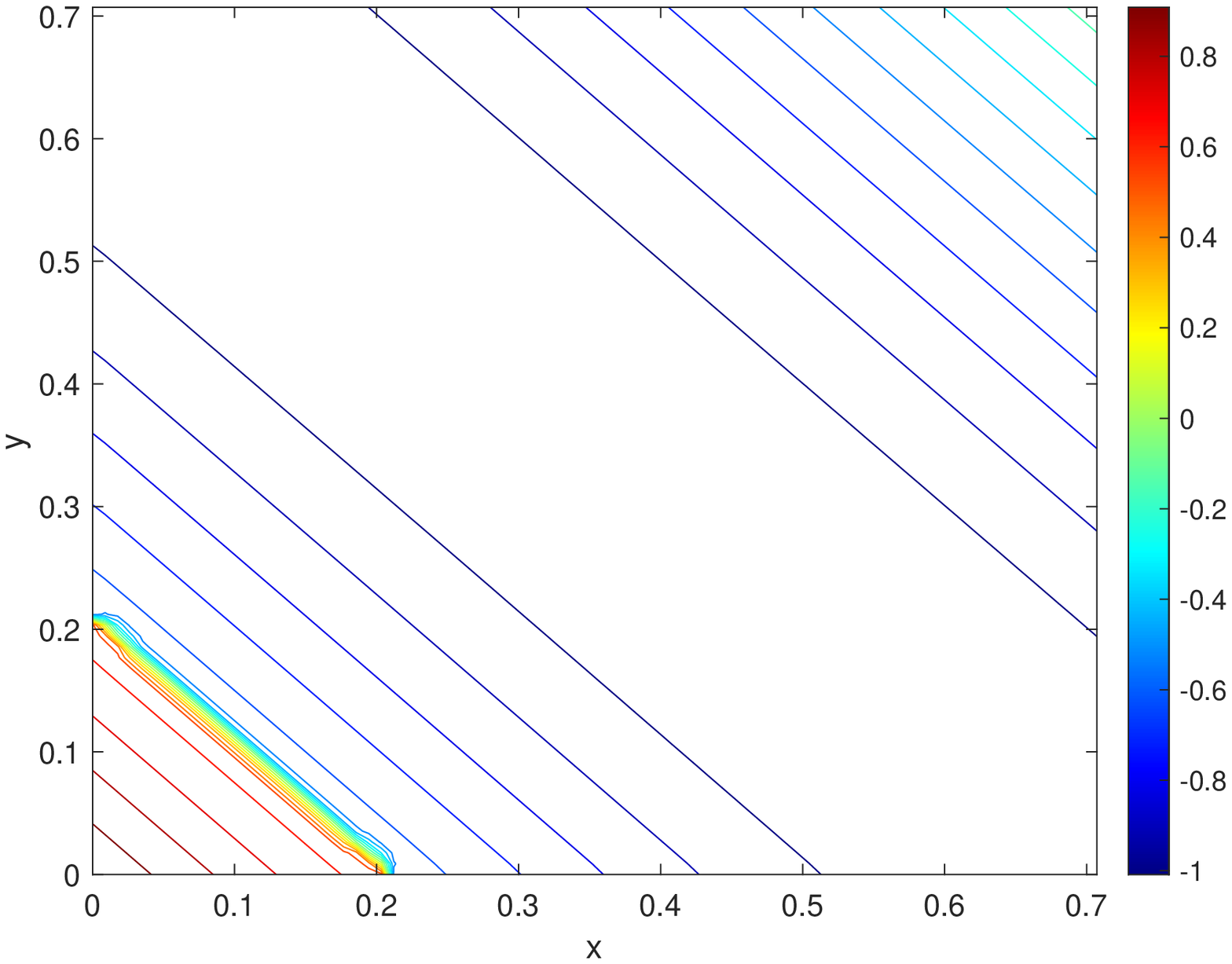}
 	\quad
    \includegraphics[width=0.45\textwidth]{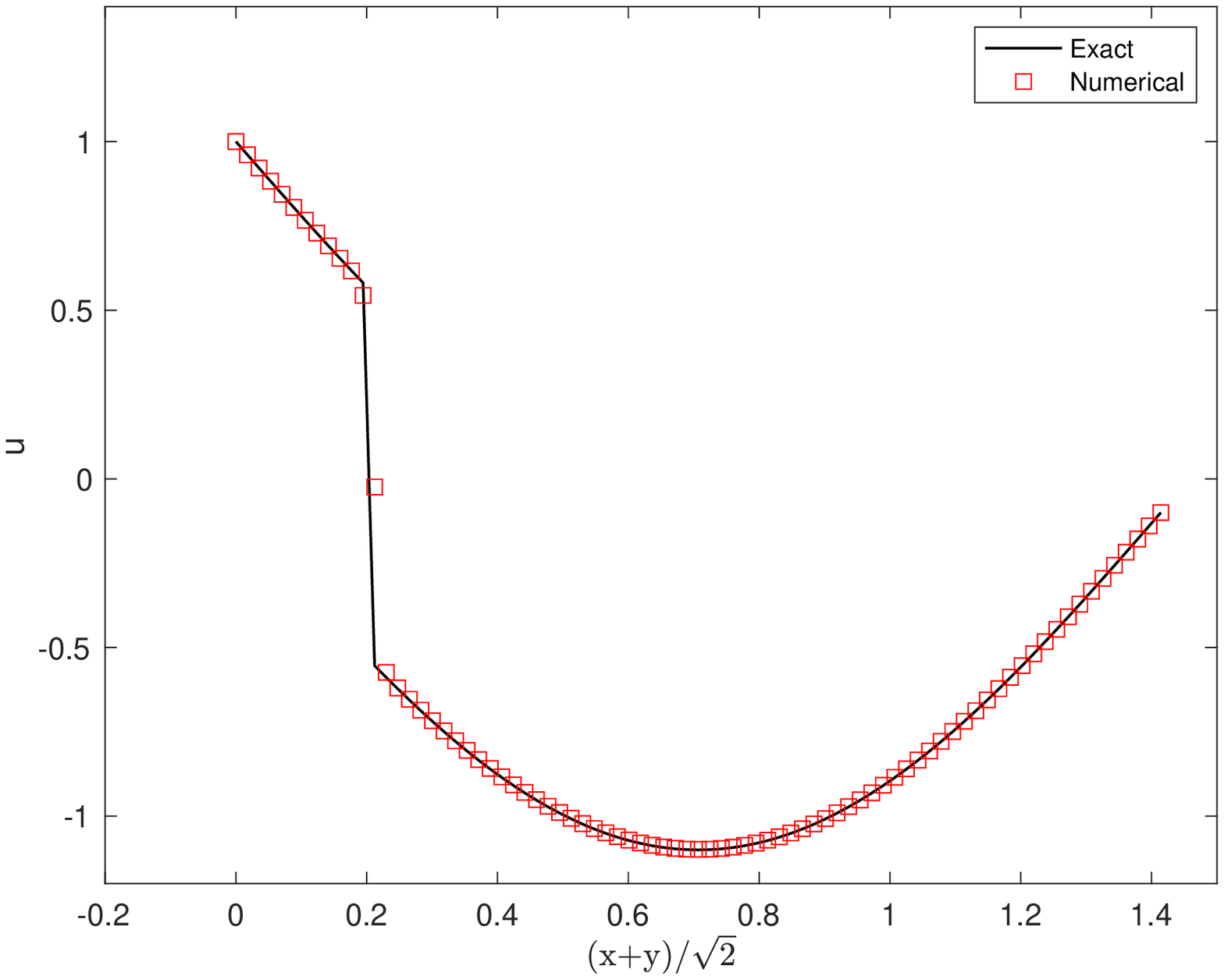}}

  \caption{{\color{orange}Example \ref{EG:UW_HWENO_BUR2D_NSM} with $80\times 80$ cells. Left: $25$ equally spaced contours of the solution from $-1.2$ to $1.1$; right: the numerical solution vs. the exact solution  along the cross-section through the northeast to southwest diagonal.}}
  \label{PIC:UW_HWENO_BUR2D_NSM}
\end{figure}
\end{exa}

\begin{exa}\label{EG:UW_HWENO_BUR2D_NSM2}
{\color{red} We consider the one-dimensional Burgers' equation viewed as a two-dimensional steady state problem}
\begin{equation}
\left(\frac{u^{2}}{2}\right)_{x} + u_{y} = 0, ~~(x, y)\in[0, 1]\times[0, 1]
\end{equation}
with the boundary conditions
\begin{equation*}
u(x, 0) = 1.5-2x, ~~u(0, y) = 1.5, ~~u(1, y) = -0.5.
\end{equation*}
{\color{orange}The exact solution consists of a fan that merges into a shock whose foot is located at $(x, y) = \left(\frac{3}{4}, \frac{1}{2}\right)$.} More precisely, the exact solution is
\[
u(x, y) = \begin{cases}
\text{if $y\geq 0.5$} & \begin{cases}
-0.5 & \text{if $-2(x - 3/4) + (y -1/2) \leq 0$}, \\
1.5 & \text{else},
\end{cases} \\
\hfil \text{else} & \text{$\max \left(-0.5, \min  \left(1.5, \frac{x-3/4}{y-1/2}\right)\right)$}.
\end{cases}
\]

This problem was studied in \cite{Cai.Gottlieb.Shu_MC1989} as a prototype example for shock boundary layer interaction. The initial condition is taken to be $u_{0}(x, y) = u_{0}(x, 0) = 1.5 - 2x$. {\color{orange}The isolines of the numerical solution and the cross-sections at $y = 0.25$ across the fan, at $y = 0. 5$ right at the junction where the fan becomes a single shock, and at $y = 0.75$ across the shock, are displayed in Figure \ref{PIC:UW_HWENO_BUR2D_NSM2}.} We can clearly observe good resolution of the numerical scheme for this example. The coefficient $\sigma$ for the dissipation \eqref{DEF:UW_HWENO2D_RESIDUAL_DISS} is taken as $2$.
\begin{figure}[!htbp]
  \centering
  \mbox{
    \includegraphics[width=0.45\textwidth]{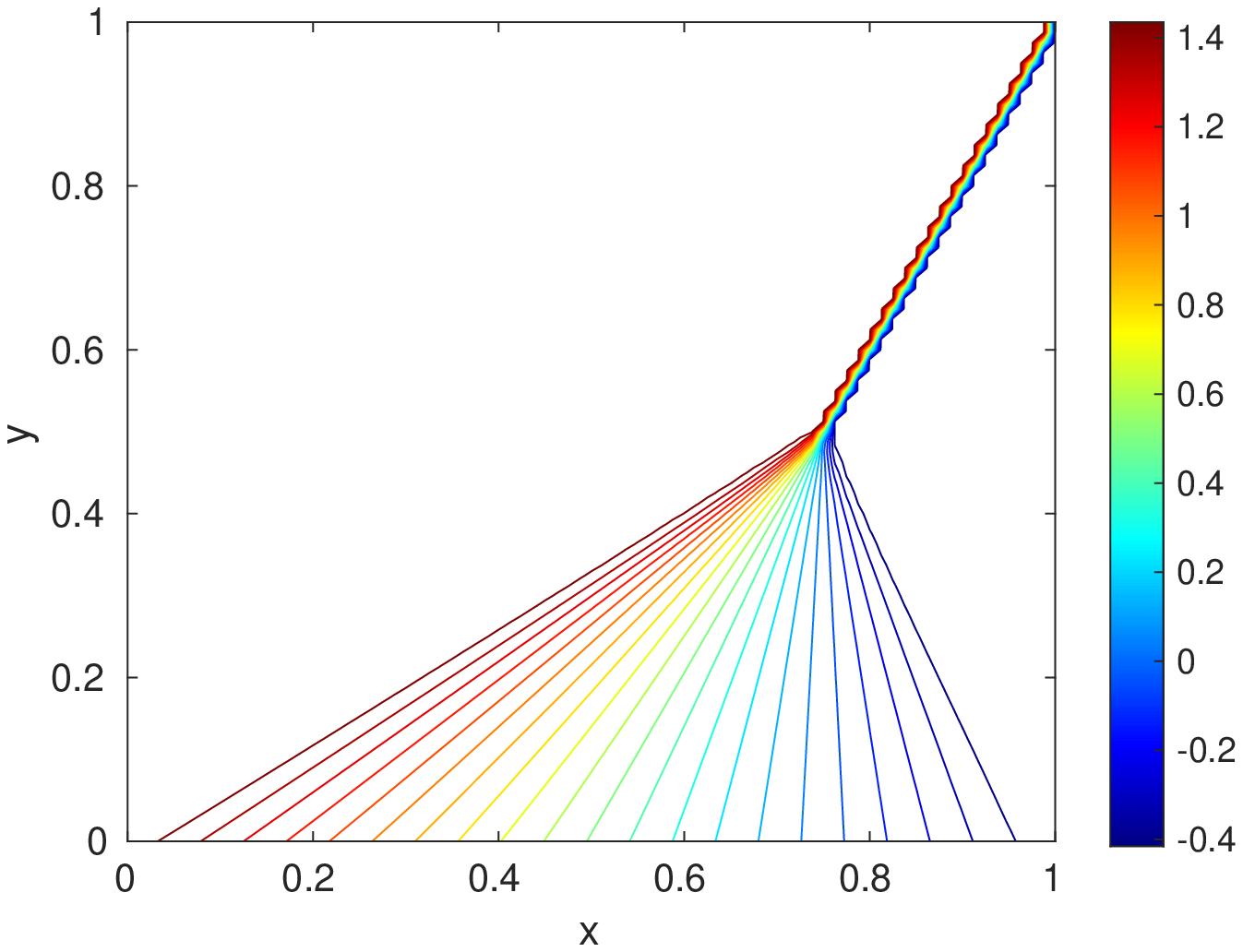}
  	\quad
    \includegraphics[width=0.45\textwidth]{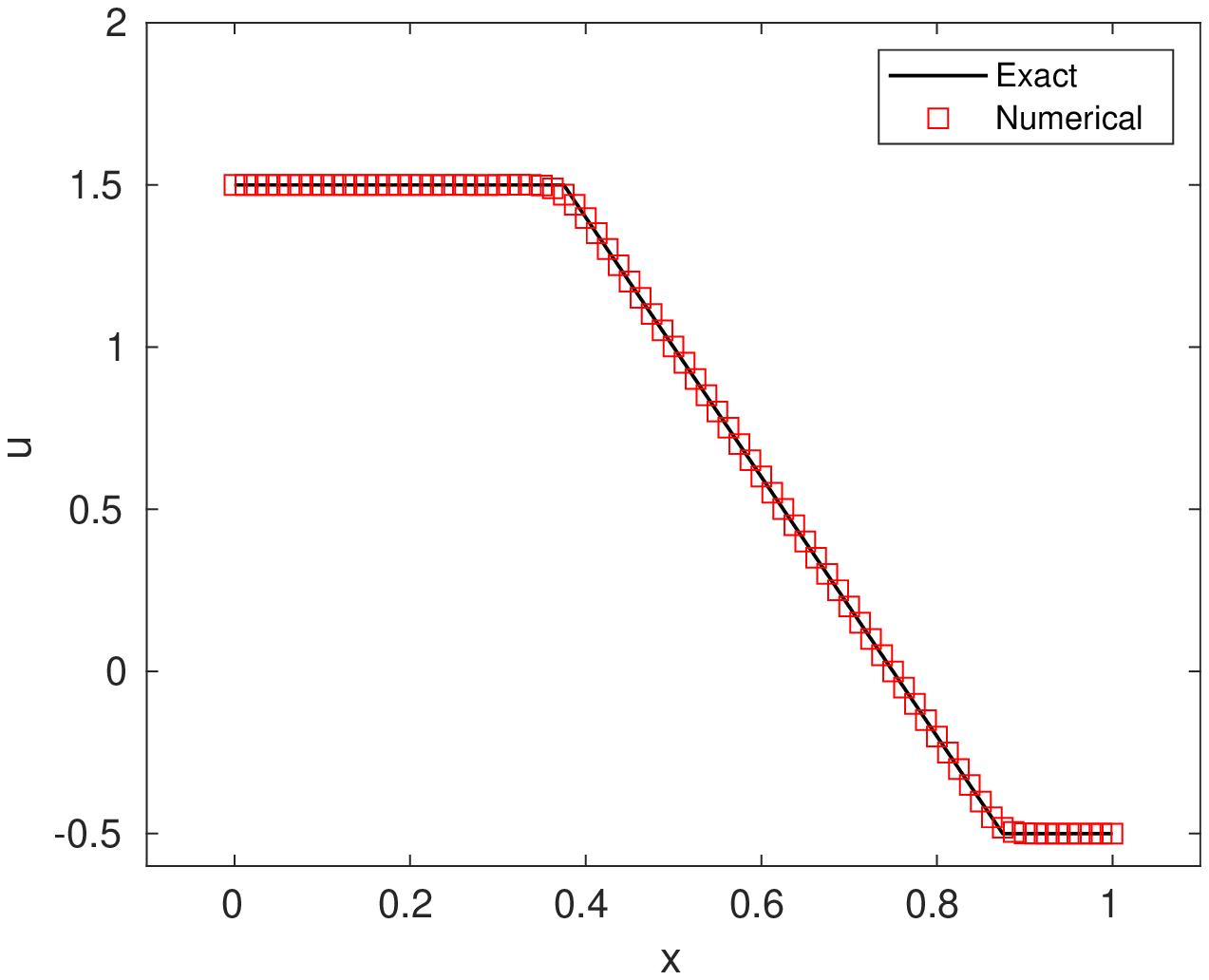}}
  \mbox{
    \includegraphics[width=0.45\textwidth]{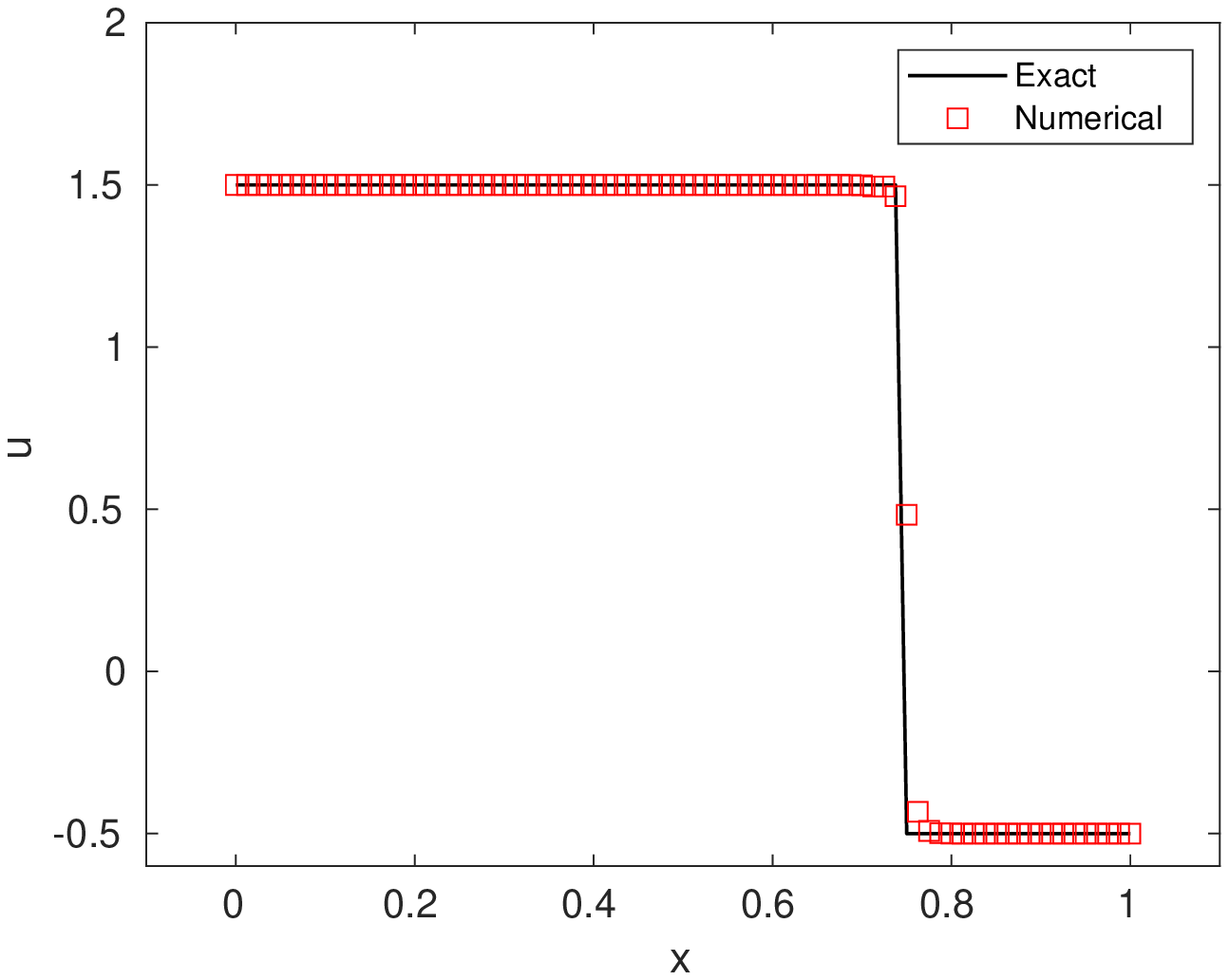}
  	\quad
    \includegraphics[width=0.45\textwidth]{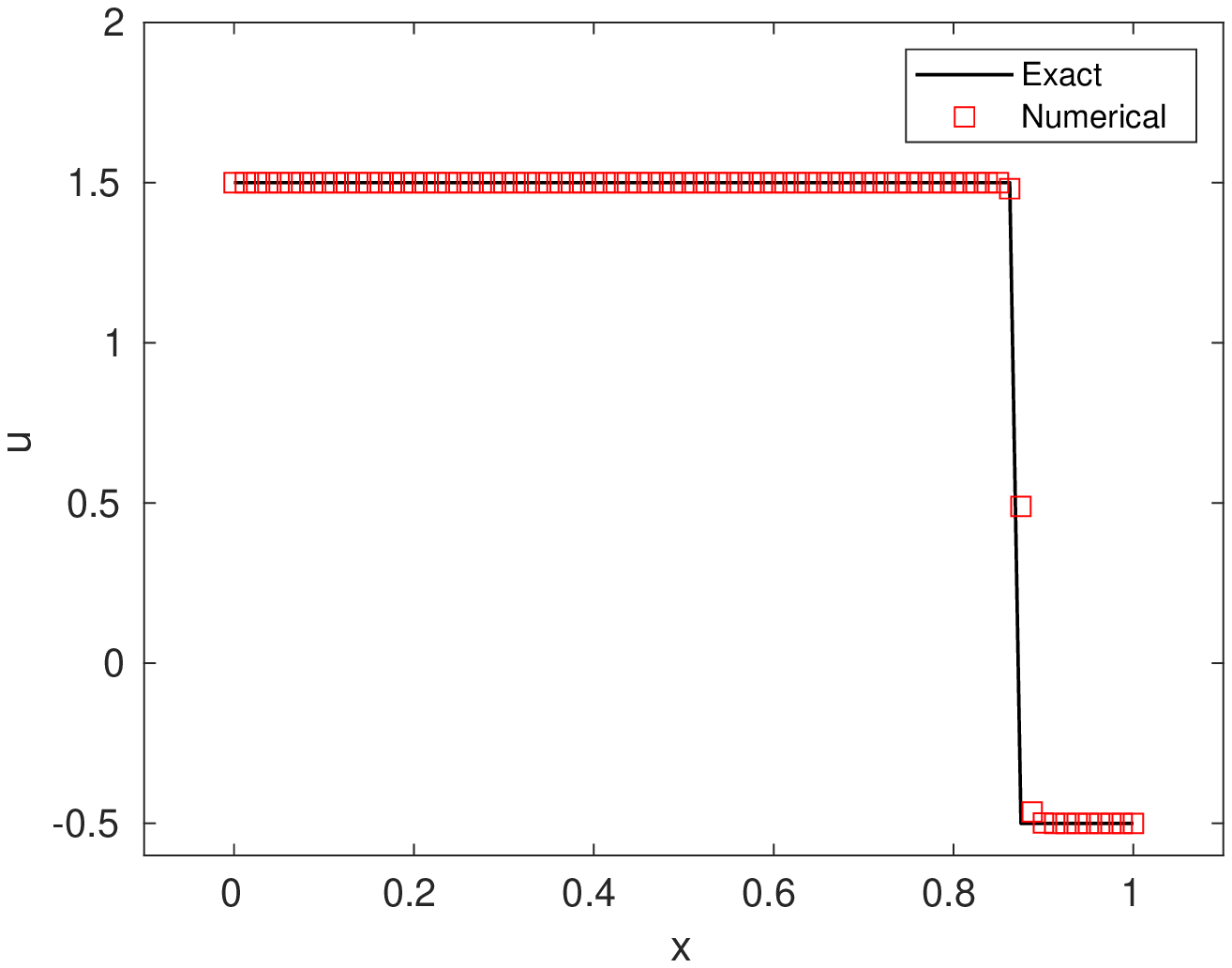}}

  \caption{{\color{orange}Example \ref{EG:UW_HWENO_BUR2D_NSM2} with $80\times 80$ cells. Top left: 25 equally spaced contour lines from -0.6 to 1.62; Top right: cross section at $y=0.25$; bottom left: cross section at $y=0.5$; bottom right: cross section at $y=0.75$.}}
  \label{PIC:UW_HWENO_BUR2D_NSM2}
\end{figure}
\end{exa}

{\color{blue}
\begin{rem}
	For the one-dimensional scalar, systems and two-dimensional scalar problems in our simulations, the affects of the parameter $\delta$ in \eqref{ROE_FUN} for the Roe's entropy correction on the problems can be ignored so that we can take a small number as $10^{-15}$. As for the choice of the parameter $\sigma$, by the definition of the dissipation residual in \eqref{DEF:UW_HWENO2D_RESIDUAL_DISS}, the magnitude of $\sigma$ decides the amount of the dissipation residual added to the vertices within each cell. If $\sigma$ is large, the numerical result will be too dissipative; otherwise, the $L^{1} $ residue is not convergent. In principle, we take a proper $\sigma$ to reach a steady state solution with a good performance of the problem.
\end{rem}
}

\begin{exa}\label{EG:UW_HWENO_CRP2D}
We consider a Cauchy-Riemann problem
\begin{equation}
\frac{\partial W}{\partial t} + A\frac{\partial W}{\partial x} + B\frac{\partial W}{\partial y} = 0, ~~ (x, y)\in [-2, 2]\times[-2, 2], ~~ t > 0,
\end{equation}
where
\begin{equation}\label{CRP_Jac}
A = \left(\begin{array}{cc}
1 & 0 \\
0 & -1
\end{array}\right)~~\text{and}~~ B = \left(\begin{array}{cc}
0 & 1 \\
1 & 0
\end{array} \right)
\end{equation}
with the following Riemann data $W = (u, v)^{T}$:
\begin{equation}\label{CRP_RMdata}
u = \begin{cases}
\hfil 1 & \text{if $x>0$ and $y>0$} \\
\hfil -1 & \text{if $x<0$ and $y>0$} \\
\hfil -1 & \text{if $x<0$ and $y<0$} \\
\hfil 1 & \text{if $x<0$ and $ y<0$}
\end{cases} ~~\text{and}~~ v = \begin{cases}
\hfil 1 & \text{if $x>0$ and $y>0$} \\
\hfil -1 & \text{if $x<0$ and $y>0$} \\
\hfil -1 & \text{if $x>0$ and $y<0$} \\
\hfil 2 & \text{if $x<0$ and $y<0$}
\end{cases}.
\end{equation}

The solution is self-similar, and therefore $W(x, y, t) = \tilde{W}\left(\frac{x}{t},\frac{y}{t}\right)$. Let $\xi = \frac{x}{t}$, $\eta = \frac{y}{t}$, then $\tilde{W}$ satisfies
\begin{equation}
(-\xi I + A)\frac{\partial\tilde{W}}{\partial\xi} + (-\eta I + B)\frac{\partial \tilde{W}}{\partial \eta} = 0,
\end{equation}
which can be written as
\begin{equation}\label{EQ:CRPs}
\frac{\partial}{\partial \xi}[(-\xi I + A)\tilde{W}] + \frac{\partial}{\partial \eta}[(-\eta I + B)\tilde{W}] =-2\tilde{W}.
\end{equation}
When $t=1$, the problem \eqref{EQ:CRPs} can be regarded as a steady state problem and solved by RD method with boundary conditions set as the exact solution and the same initial condition as in \eqref{CRP_exact}. The coefficient $\sigma$ for the dissipation \eqref{DEF:UW_HWENO2D_RESIDUAL_DISS} is taken as 1 and the parameter $\delta $ in \eqref{ROE_FUN} for the Roe's entropy correction is taken as  0.4. {\color{orange}The numerical results are shown in Figure \ref{PIC:UW_HWENO_CRP2D} and the convergence history of $L^{1}$ residue stagnates around $10^{-7}$ level, as shown in Figure \ref{PIC:UW_HWENO_CRP2D_LOG10RES}.}
\begin{align}\label{CRP_exact}
u = \begin{cases}
\hfil 1 & \text{if $x>1$ and $y>1$} \\
\hfil -1 & \text{if $x>1$ and $y<1$} \\
\hfil -1 & \text{if $x<1$ and $y>1$} \\
\hfil 1.5 & \text{if $x<1$ and $-1<y<1$}\\
\hfil 1 & \text{if $x<1$ and $y<-1$}
\end{cases} ~~\text{and}~~ v = \begin{cases}\hfil
\hfil 1 & \text{if $x>-1$ and $y>1$}\\
\hfil -1 & \text{if $x<-1$ and $y>1$} \\
\hfil -1 & \text{if $x>-1$ and $y<1$} \\
\hfil 1.5 & \text{if $x<-1$ and $-1< y <1$} \\
\hfil 2 & \text{if $x<-1$ and $y<-1$}
\end{cases}.
\end{align}

\begin{figure}[!htbp]
  \centering
  \mbox{
    \includegraphics[width=0.45\textwidth]{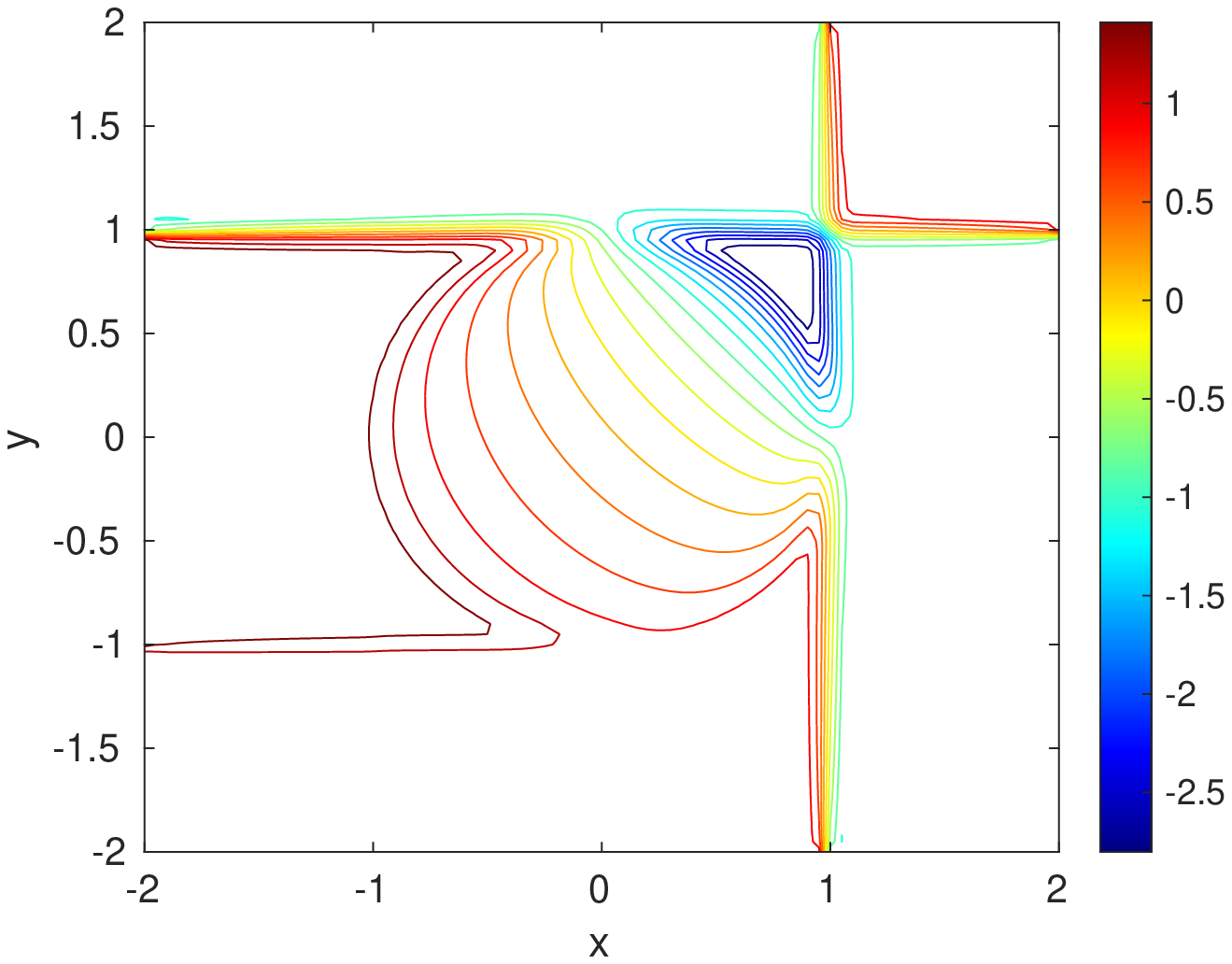}
  	\quad
    \includegraphics[width=0.45\textwidth]{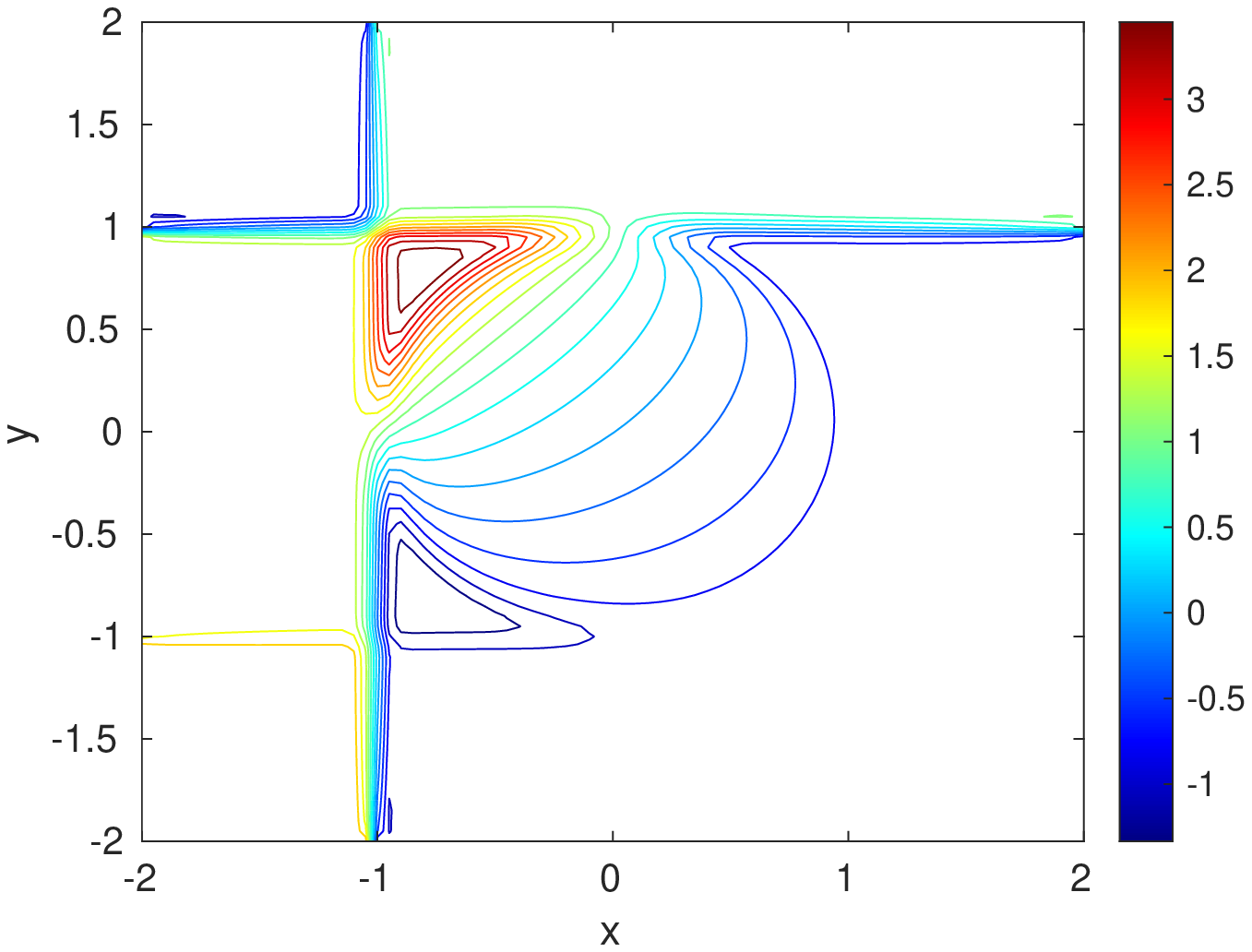}}

  \caption{{\color{orange}Example \ref{EG:UW_HWENO_CRP2D} with $80 \times 80$ cells. Left: $20$ equally spaced contours for $u$ from $-3.05$ to $1.66$; right: $20$ equally spaced contour for $v$ from $-1.6$ to $3.45$.}}
  \label{PIC:UW_HWENO_CRP2D}
\end{figure}

\begin{figure}[!htbp]
	\centering
	\includegraphics[width=0.45\textwidth]{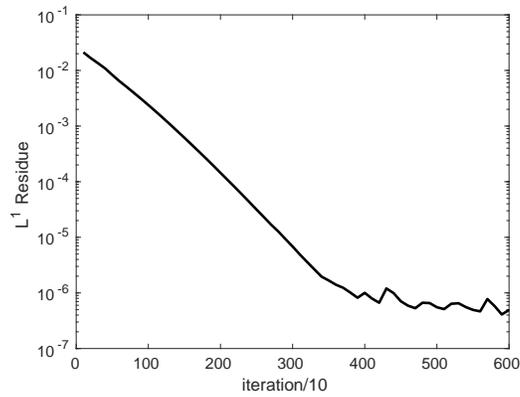}	
	\caption{{\color{orange}The convergence history of $L^{1}$ residue for the Cauchy Riemann problem in Example \ref{EG:UW_HWENO_CRP2D}.}}
	\label{PIC:UW_HWENO_CRP2D_LOG10RES}
\end{figure}
\end{exa}

\begin{exa}\label{EG:UW_HWENO_SHOCKREF2D}
We consider a regular shock reflection problem of the steady state solution of the two-dimensional Euler equations
\begin{equation}
{\bf f}({\bf u})_{x} + {\bf g}({\bf u})_{y} = 0, ~~(x, y)\in [0, 4]\times[0, 1],
\end{equation}
where ${\bf u} = (\rho, \rho u, \rho v, E)^{T}$, ${\bf f}({\bf u}) = (\rho u, \rho u^{2} + p, \rho uv, u(E + p))^{T}$, and ${\bf g}({\bf u}) = (\rho v, \rho uv, \rho v^{2} + p, v(E + p))^{T}$. Here $\rho$ is the density, $(u, v)$ is the velocity, $E$ is the total energy and $p = (\gamma -1)(E - \frac{1}{2}(\rho u^{2} + \rho v^{2}))$ is the pressure. $\gamma$ is the gas constant which is again taken as $1.4$ in our numerical tests.

The initial condition is taken to be
\[
(\rho, u, v, p) = \begin{cases}
(1.69997, 2.61934, -0.50632, 1.52819) & \text{on $y = 1$},\\
\hfil (1, 2.9, 0, \frac{1}{\gamma}) & \text{otherwise}.
\end{cases}
\]
The boundary conditions are given by
{\color{orange}
\[
(\rho, u, v, p) = (1.69997, 2.61934, -0.50632, 1.52819) ~~ \text{on $y = 1$},
\]}
and reflective boundary condition on $y = 0$. The left boundary at $x = 0$ is set as inflow with $(\rho, u, v, p)=(1, 2.9, 0, \frac{1}{\gamma})$, and the right boundary at $x = 4$ is set to be an outflow with no boundary conditions prescribed. {\color{orange}The numerical results are shown in Figure \ref{PIC:UW_HWENO_SHOCKREF2D}.} We can clearly see the good resolutions of the incident and reflected shocks.  {\color{orange}The coefficient $\sigma$ for the dissipation \eqref{DEF:UW_HWENO2D_RESIDUAL_DISS} is taken as 8. The convergence history of $L^{1}$ residue is shown in Figure \ref{PIC:UW_HWENO_SHOCKREF2D_LOG10RES}.}
\begin{figure}[!htbp]
	\centering
	\mbox{
		\includegraphics[width=0.45\textwidth]{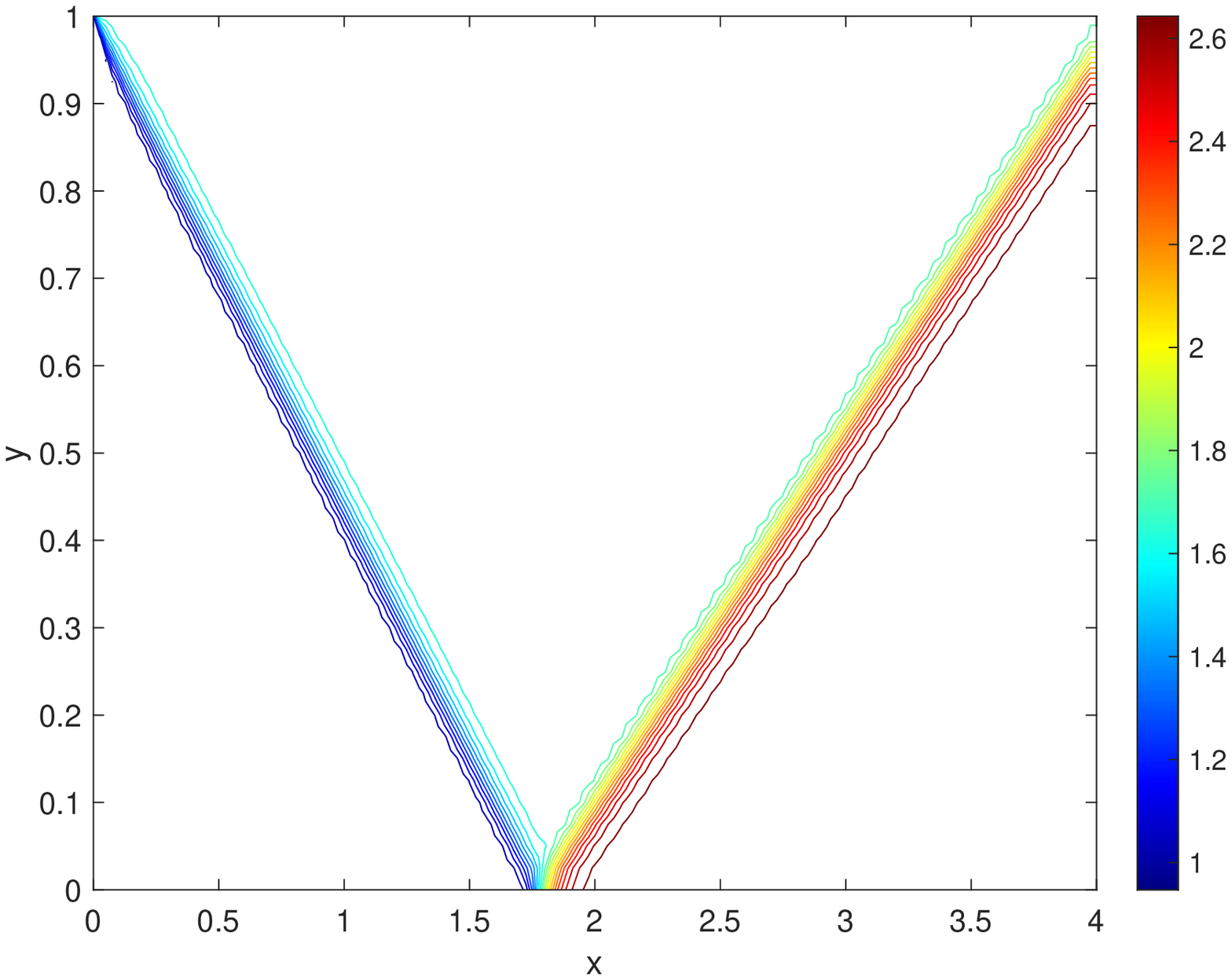}
		\quad
		\includegraphics[width=0.45\textwidth]{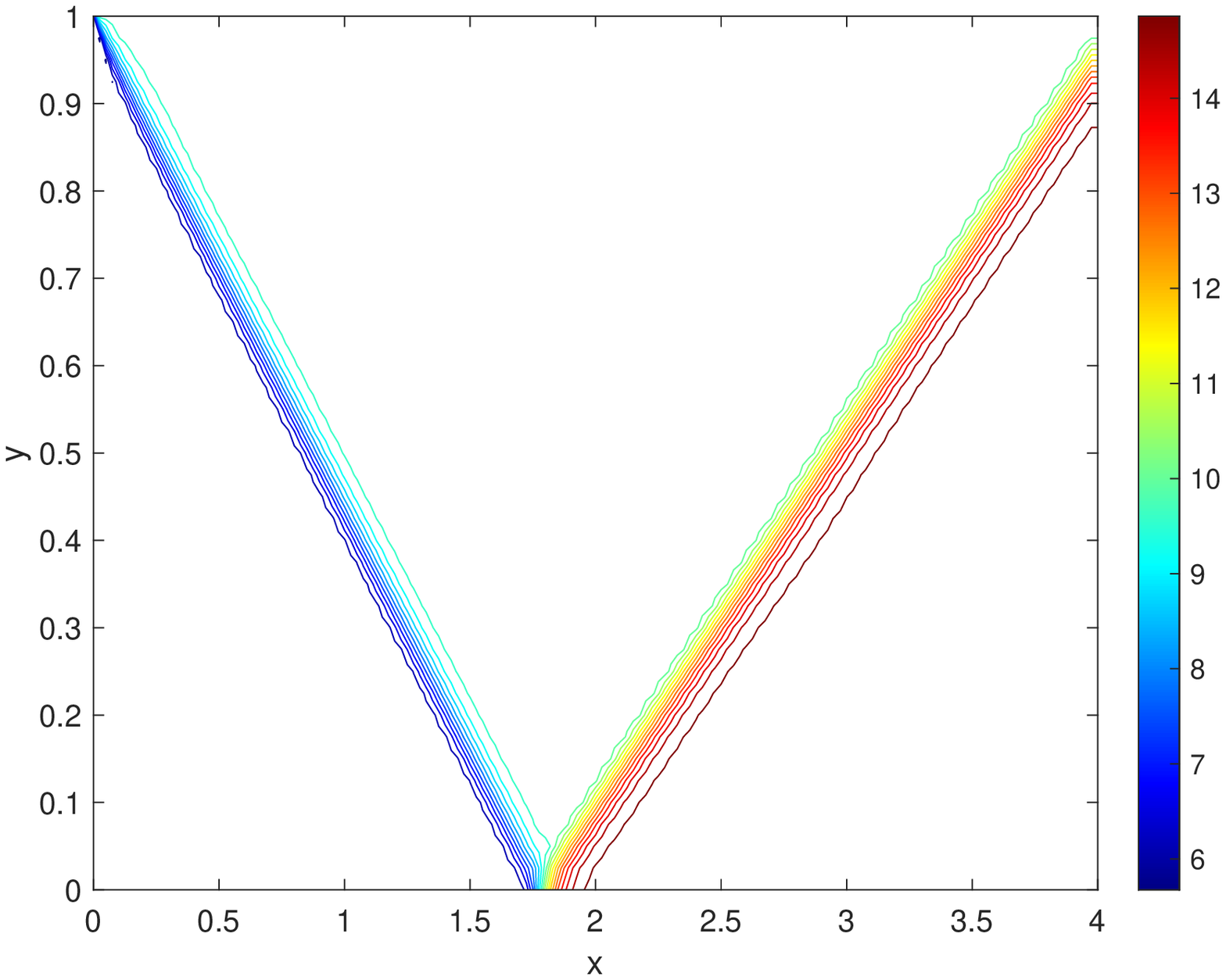}}
	
	\caption{{\color{orange}Example \ref{EG:UW_HWENO_SHOCKREF2D} with $160 \times 40$ cells. Left: $25$ equally spaced contours for the density from $0.87$ to $2.72$; right: $25$ equally spaced contour for the energy from $4.8$ to $15.3$.}}
	\label{PIC:UW_HWENO_SHOCKREF2D}
\end{figure}

\begin{figure}[!htbp]
	\centering
	\includegraphics[width=0.45\textwidth]{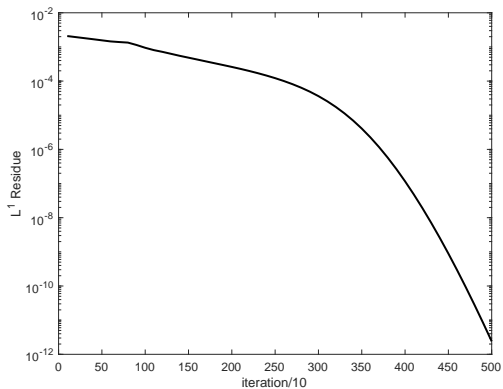}	
	\caption{{\color{orange}The convergence history of $L^{1}$ residue for the shock reflection problem in Example \ref{EG:UW_HWENO_SHOCKREF2D}.}}
	\label{PIC:UW_HWENO_SHOCKREF2D_LOG10RES}
\end{figure}
\end{exa}

{\color{blue}
\begin{rem}
	In two-dimensional system cases, the magnitude of the parameter $\delta$ in \eqref{ROE_FUN} for the Roe's entropy correction affects the capability of preserving the non-oscillatory property. In our simulations, we choose a proper $\delta$ to avoid the spurious oscillations and also $\sigma$ which makes the numerical result not too dissipative.
\end{rem}
}

\section{Conclusion}
{\color{orange}A high order RD  conservative finite difference HWENO method was proposed for solving steady state hyperbolic equations with source terms on uniform meshes. The method is based on a novel HWENO  scheme to achieve high order accuracy.} {\color{orange}Compared with WENO scheme, the advantage of traditional HWENO is that the stencil is more compact and the error is smaller under the same grid and accuracy. However, the additional auxiliary equations are required in traditional HWENO framework, and it is not clear how to distribute residuals for the auxiliary equations when solving two-dimensional steady state conservation laws.  In this paper, the novel HWENO framework developed in \cite{renhigh} is extended to the steady state hyperbolic conservation law, and the residual distribution method is developed. The new framework not only inherits the advantages of the traditional HWENO, but also does not need to introduce any auxiliary equations, which leads to less storage and a low-cost advantage.} {\color{orange}We apply this proposed method to both scalar and system problems in one and two dimensions including  Burgers' equation}, shallow water equations, nozzle flow problem, Cauchy Riemann problem and Euler equations. In all simulations, we observe that we get the sixth order in smooth regions, respectively, and clearly see high resolutions around a shock. The extension to unsteady problems will be explored in future.

\section*{Acknowledgments}
The research was started in 2016 when J. Lin was a visiting Ph.D. student at University of Z\"urich in Switzerland and she was also partly supported by SNF FZEB-0-166980 grant and National Science Foundation (China) grant 12071392. And thank Professor Xinghui Zhong at Zhejiang University for her valuable comments on the paper.

\appendix

\bibliographystyle{abbrv}
\bibliography{reference/myrefer}
\end{document}